\documentclass[12pt]{amsart}

\usepackage{amssymb,euscript,amsmath, mathrsfs,amscd}
\usepackage{epsfig} 
\usepackage{enumerate}
\usepackage{fullpage}
\usepackage{tikz}

\newcounter{ENUM}
\newcommand{\beas}{\begin{eqnarray*}}
\newcommand{\eeas}{\end{eqnarray*}}

\newtheorem{thm}{Theorem}[section]
\newtheorem{prop}[thm]{Proposition}
\newtheorem{lem}[thm]{Lemma}

\newtheorem{cor}[thm]{Corollary}
\newtheorem{conj}[thm]{Conjecture}

\theoremstyle{definition}
\newtheorem{defn}[thm]{Definition}

\newtheorem{ex}[thm]{Example}

\newtheorem{notn}[thm]{Notation}
\newtheorem{rem}[thm]{Remark}

\numberwithin{equation}{section}

\def\Q{{\mathbb Q}}

\def\R{{\mathbb R}}
\def\fS{{\mathfrak S}}

\newcommand{\bm}[1]{{\boldsymbol{#1}}}
\def\0{\bm{0}}

\def\c{\bm{c}}
\def\e{\bm{e}}
\def\k{\bm{k}}
\def\m{\bm{m}}

\def\v{\bm{v}}
\def\u{\bm{u}}

\def\l{\ell}

\newcommand{\wt}{\mathrm{wt}}
\newcommand{\intwt}{\mathrm{intwt}}
\newcommand{\bdwt}{\mathrm{bdwt}}
\newcommand{\Lat}{\mathrm{Lat}}

\newcommand{\prune}{\mathrm{Prune}}
\newcommand{\nvol}{\mathrm{nvol}}

\def\ncone{\operatorname{ncone}}
\def\conv{\operatorname{ConvexHull}}
\def\cone{\operatorname{Cone}}

\def\Perm{\operatorname{Perm}}

\def\td{\operatorname{Td}}

\renewcommand{\k}{\textbf{k}}

\newcommand{\B}{{\mathcal B}}
\newcommand{\D}{{\Gamma}}
\newcommand{\spd}{\mathcal F}
\newcommand{\lspd}{\mathcal O}
\newcommand{\Sp}{\mathrm{Sp}}

\newcommand{\Bi}{\mathrm{Binom}}
\newcommand{\sgn}{\mathrm{sgn}}

\newcommand{\Tdcoeff}{\mathrm{Tdcoeff}}

\newcommand{\hB}{\overline{\mathcal B}}
\newcommand{\CC}{\mathcal{C}}

\newcommand{\bv}{\alpha^\textsc{bv}}

\newcommand{\NN}{{\mathbb{N}}}

\newcommand{\T}{\mathcal{T}}
\newcommand{\ZZ}{\mathbb{Z}}

\newcommand\commentout[1]{}
\author{Federico Castillo}
\address{Federico Castillo, Department of Mathematics, University of
  Kansas, 1450 Jayhawk Blvd, Lawrence, KS 66045 USA.}
  \email{fcastillo@ku.edu}
  
  \author{Fu Liu}
  \address{Fu Liu, Department of Mathematics, University of California, Davis, One Shields Avenue, Davis, CA 95616 USA.}
\email{fuliu@math.ucdavis.edu}

\begin{document}
\title{On the Todd Class of the Permutohedral variety}

\maketitle
\begin{abstract}
In the special case of braid fans, we give a combinatorial formula for the Berline-Vergne's construction for an Euler-Maclaurin type formula that computes number of lattice points in polytopes. Our formula is obtained by computing a symmetric expression for the Todd class of the permutohedral variety.
By showing that this formula does not always have positive values, we prove that the Todd class of the permutohedral variety $X_d$ is not effective for $d\geq 24$.

Additionally, we prove that the linear coefficient in the Ehrhart polynomial of any lattice generalized permutohedron is positive. 
\end{abstract}

\section{Introduction}
\label{sec:in}

Let $\Lambda$ be a lattice of finite rank and $V=\Lambda\otimes\R$ be the corresponding real finite-dimensional vector space. A \textbf{lattice polyope} in $V$ is a polytope such that all of its vertices lie in $\Lambda$. 
A classical problem in the crossroads between enumerative combinatorics and discrete geometry is that of counting lattice points in lattice polytopes. For any polytope $P\subset V$ we define $\Lat(P):=|P\cap\Lambda|$. One of the earliest results in the area is Pick's theorem, which says that for any lattice polygon $P\subset\R^2$ we have
\[
\Lat(P) = a(P)+\dfrac{1}{2}b(P)+1,
\]
where $a(P)$ is the area of $P$ and $b(P)$ is the number of lattice points on the boundary of $P$. One way to obtain a higher dimensional analog of Pick's formula is to find a formula relating the number of lattice points of $P$ with the different normalized volumes of the faces $F$ of $P$. We want a real-valued function $\alpha$ on pairs $(F,P)$, where $F$ is a face of a lattice polytope $P$, such that
\begin{equation}\label{eq:exterior}
	\Lat(P) = \sum_{F: \text{ a face of $P$}} \alpha(F,P) \ \nvol(F),
\end{equation}
where $\nvol(F)$ is the normalized volume of $F$.
It is clear that for a given lattice polytope $P$ one can always find many functions $\alpha$ satisfying \eqref{eq:exterior}. What we want is a function that works simultaneously for \emph{all} lattice polyopes. We can do this by requiring the function $\alpha$ to be \emph{local}, i.e., the value of $\alpha(F,P)$ only depends on the local geometry of $P$ around $F$, or more specifically, the value only depends on $\ncone(F,P)$, the normal cone of $P$ at $F.$ Any \emph{local} function $\alpha$ that satisfies Equation \eqref{eq:exterior} for all lattice polytopes $P$ is called a \textbf{McMullen function}, since McMullen was the first to prove their existence \cite{mcmullenext}. His proof is nonconstructive and shows that there are infinitely many McMullen functions.
In the present paper we compute the values for a particular McMullen function on a special family of polytopes: generalized permutohedra, originally defined by Postnikov as deformations of usual permutohedra. An alternative and equivalent definition can be given in terms of normal fans: The \textbf{braid fan} $\Sigma_d$ is the complete fan in the quotient space $W_d := \R^{d+1}/(1,1,\dots,1)$ obtained from the hyperplane arrangement $H_{i,j}:=\{x\in\R^{d+1}: x_i-x_j=0\}$ for $1\leq i<j\leq d+1$.
A \textbf{generalized permutohedron} is a polytope whose normal fan is a coarsening of the braid fan $\Sigma_d.$

 Our methods for computing a McMullen function for generalized permutohedra are based on the theory of toric varieties.
\subsection{Todd classes of toric varieties}\label{subsec:toddintro}

Let $P$ be a lattice polytope with normal fan $\Sigma$ and $X_\Sigma$ be the associated toric variety. The Todd class $\td(X_\Sigma)$ is an element in the Chow ring of $X_\Sigma$. As such it can be written as a $\Q$-linear combination of the toric invariant cycles $[V(\sigma)]$:
\begin{equation}\label{eq:toddexpression}
\td(X_\Sigma) = \sum_{\sigma\in \Sigma} r_\Sigma(\sigma)\ [V(\sigma)],\quad r_\Sigma(\sigma)\in\Q.
\end{equation}
Since the cycles $[V(\sigma)]$ satisfy algebraic relations, the values $r_\Sigma(\sigma)$ satisfying \eqref{eq:toddexpression} are not uniquely determined. 
An amazing connection with lattice polytopes is given by the fact that any function $r_\Sigma(\cdot)$ satisfying \eqref{eq:toddexpression} defines a function $\alpha$ satisfying \eqref{eq:exterior} for $P$ by setting 
\[\alpha(F,P)=r_{\Sigma}(\ncone(F,P)).\] A proof of this fact can be found in Danilov's 1978 survey \cite{danilov} where he further asked if there exist a function $r$ that depends only on the cone $\sigma$ and not on $\Sigma$, in other words, if there exist a \emph{local} function $r$ satisfying Equation \eqref{eq:toddexpression} for all fans $\Sigma$.
Accordingly, we call such a function $r$ on pointed cones a \emph{Danilov function}. By setting \[\alpha(F,P)=r(\ncone(F,P)),\] any Danilov function gives a McMullen function.

We want to briefly remark on two constructions of Danilov functions from the last two decades. Pommersheim and Thomas \cite{toddclass} gave a construction of a Danilov function $r(\sigma)$ that depends on choosing a \emph{complement map} for subspaces. Originally they do this by choosing a complete flag of subspaces, which has a technical issue that their construction of $r(\cdot)$ is only defined for cones that are ``generic'' with respect to the chosen flag. Hence strictly speaking, their function $r$ is only an ``almost'' Danilov function.

A couple of years later Berline and Vergne \cite{bvoriginal} constructed a McMullen function with the property that it is computable in polynomial time fixing the dimension and it is a valuation on cones. We call this construction the BV-function, and denote it by $\bv$. 

Later in \cite{bvtodd}, 
they showed that if a function $r$ on pointed cones is defined by
\[ r(\sigma) = \bv\left( F,P \right) \text{ as long as $\sigma = \ncone(F,P)$,}\]
then it is a Danilov function.
For convenience, we abuse the notation, and consider $\bv$ to be both a function on pairs $(F,P)$ and a function on cones with the connection that
\[ \bv(F,P) = \bv(\ncone(F,P)).\] 
Thus, $\bv$ is both a McMullen function and a Danilov function.
In \cite{garpomm} Pommersheim and Garoufalidis proved that using an inner product for a complement map in the methods of \cite{toddclass} results in the Danilov function $\bv$, which in turns gives an alternative way of computing it. 

Both constructions, Berline-Vergne's and Pommersheim-Thomas', are algorithmic. A priori it is very hard to get formulas for general cones. There are very few examples of fans $\Sigma$ for which $\bv(\sigma)$ (or any other Danilov function) have been computed for all $\sigma\in\Sigma$. In this paper, we focus on computing the BV-function on all cones in braid fans using tools developed in previous work of the authors. (See Section \ref{subsec:braid} for the definition of braid fans.)

In \cite{bvalpha} we exploited an extra symmetry property satisfied by the function $\bv$, and used this symmetry to study the values on cones in braid fans. One main result in \cite{bvalpha} is the \emph{uniqueness theorem}, which in the context of the present paper states that, for the specific example of braid fans, $\bv$ is the unique function satisfying Equation \eqref{eq:toddexpression} and being invariant under the permutation action of the symmetric group on the ambient space.
Using this, we obtain the main result of this paper - Theorem \ref{thm:main} - which gives a combinatorial formula for $\bv$ on all cones in braid fans.

\subsection{Connection to Ehrhart theory}
In \cite{ehrhart} Ehrhart proved that for every lattice polytope $P$ the function $\Lat(tP)$ for $t\in \NN$ is a polynomial in $t$ of dimension $d=\dim P.$ More precisely, there exist $a_0, a_1, \dots, a_d \in\mathbb{Q}$ such that for all $t \in \NN$
\[
\Lat(tP) = a_0+a_1t^1+a_2t^2+\cdots+a_dt^d.
\]

The right hand side is called the \textbf{Ehrhart polynomia}l of $P$. Given a McMullen formula $\alpha$ one can deduce that 
\begin{equation}\label{eq:refine}
	a_k = \sum_{\substack{F: \text{a face of $P$} \\\dim F=k}}\alpha(F,P) \ \nvol(F). 
\end{equation}

The first, second, and the last coefficients in the Ehrhart polynomial of a lattice polytope are well-understood. In particular, they are all positive. However, any of the remaining ``middle coefficients'' $a_{d-2}, a_{d-3}, \dots, a_1$ can be negative. 
We call a lattice polytope $P$ \textbf{Ehrhart positive} if all the (middle) coefficients of its Ehrhart polynomial are positive (see \cite{positivesurvey} for a recent survey on Ehrhart positivity). 
One of the main motivations for \cite{bvalpha} was to prove a conjecture of De Loera et al. asserting that matriod polytopes are Ehrhart positive \cite{deloera}. 
Noticing that matroid polytopes belong to the family of generalized permutohedra, we focus on the latter larger family of polytopes.
\begin{conj}[Conjecture 1.2 of \cite{bvalpha}]\label{conj:main}
Lattice generalized permutohedra are Ehrhart positive.
\end{conj}

One observes that a consequence of Equation \eqref{eq:refine} is that if we have a McMullen function $\alpha$ such that $\alpha(F,P)$ is positive for all faces $F\subset P$ then $P$ is Ehrhart positive. (The converse is not true as shown in Section 3.4 of \cite{nill}.) Using the fact that the BV-function $\bv$ is a McMullen function and it has certain valuation properties, we showed in \cite{bvalpha} that the following conjecture (if true) implies Conjecture \ref{conj:main}.
\begin{conj}[Conjecture 1.3 of \cite{bvalpha}]\label{conj:refined}
Let $P$ be a generalized permutohedron and $F\subset P$ a face, then $\bv(F,P)$ is positive. 
Equivalently, $\bv(\sigma)$ is positive for every cone $\sigma$ in the braid fan.
\end{conj}

In \cite{bvalpha}, we were able to prove that the third and fourth coefficients of the Ehrhart polynomial of any lattice generalized permutohedron are positive by showing that $\bv(F,P)$ is positive for any pair $(F,P)$ in which $P$ is a lattice generalized permutohedron and $F$ is a face of $P$ of codimension at most $3,$ which is equivalent to that $\bv(\sigma)$ is positive for any cone $\sigma$ of dimension at most $3$ in a braid fan. 

Despite of these positive results we've obtained in our previous work towards Conjecture \ref{conj:refined}, in the present paper we use our main result - the combinatorial formula described in Theorem \ref{thm:main} - to find negative values for $\bv$ on some cones in braid fans, hence disproving Conjecture \ref{conj:refined}. Note that this does not imply that Conjecture \ref{conj:main} is false, and in fact we present a proof, independent of the rest of the paper, that the linear coefficient of the Ehrhart polynomial of any lattice generalized permutohedron is positive, providing further evidence to Conjecture \ref{conj:main}. This positivity result of linear Ehrhart coefficient was proved independently by Jochemko and Ravichandran in \cite{kat} using different techniques from ours. 
More recent evidence for Conjecture \ref{conj:main} was found recently by Ferroni in \cite{hypersimplices} where it is proved that hypersimplices are Ehrhart positive. So even though we disproved our Conjecture \ref{conj:refined} we still believe Conjecture \ref{conj:main} is true, but its resolution requires a different approach.

Finally, as a consequence of these negative $\bv$-values we also obtained the following result about the permutohedral variety of independent interest.

\begin{thm}\label{thm:noneffective}
The Todd class of the permutohedral variety $X_d$ is not effective for $d\geq 24$. That is, there is no way of expressing it as a nonnegative combination of cycles.
\end{thm}

\subsection{Organization}
This paper is organized as follows. In Section \ref{sec:prelim} we set the preliminaries about toric varieties in an elementary way.
In Section \ref{sec:spider} we define combinatorial diagrams that will be used to express formulas asserted in our main result Theorem \ref{thm:main}. In Section \ref{sec:gen} we do some computations in the Chow ring that lead to our explicit general formula in Section \ref{sec:bv}. Section \ref{sec:ex} contains applications of our main theorem. Finally in Section \ref{sec:edge} we prove that the linear coefficient in the Ehrhart polynomial of any lattice generalized polyhedron is positive. 

\section*{Acknowledgements}
The second author is partially supported by a grant from the Simons Foundation \#426756. This project started when both authors were attending the program ``Geometric and Topological Combinatorics'' at the Mathematical Sciences Research Institute in Berkeley, California, during the Fall 2017 semester, and they were partially supported by the NSF grant DMS-1440140.

We thank Jamie Pommersheim for many helpful discussions, and thank Jos\'e Gonz\'alez for enlightening conversations about toric varieties.

\section{Preliminaries and notation.}\label{sec:prelim}
We assume familiarity with the concepts of polytopes, normal fans, and toric varieties. A good reference is \cite{ewald}.  Here we review concepts and notation that we are going to use. As standard we denote $[d+1]:=\{1,2,3,\cdots,d,d+1\}$. The set of all subsets of $[d+1]$ form a poset $\B_{d+1}$ called the \textbf{boolean algebra} and we define the \textbf{truncated boolean algebra}, denoted by $\hB_{d+1}$, to be the poset obtained from $\B_{d+1}$ by removing $[d+1]$ and $\emptyset$. Two elements $T,T'\in \hB_{d+1}$ are \textbf{incomparable} if neither $T\subseteq T'$ nor $T\supseteq T'$.

\begin{notn}
From here on we use the symbol $\subset$ to denote \emph{proper} subset, instead of $\subsetneq$.
\end{notn}

A \textbf{$k$-chain} $T_{\bullet}=(T_1,\cdots,T_k)$ is a sequence of $k$ totally ordered elements of $\hB_{d+1}$.
The set of all $k$-chains in $\hB_{d+1}$ is denoted $\CC_{d+1}^k$ and let $\CC_{d+1}=\bigcup_k \CC_{d+1}^k$.

\subsection{Braid fan and Permutohedral variety}\label{subsec:braid}
Let $V_d$ be the $d$-dimensional real vector space $\textbf{1}^\perp\subset \R^{d+1}$, where $\textbf{1}$ is the all one vector. Its dual is $W_d = \R^{[d+1]}/(\textbf{1})$. 

\begin{defn}
  Given a point $\v = (v_1,v_2,\cdots,v_{d+1}) \in \R^{d+1}$, we define the \textbf{usual permutohedron}
\[\Perm(\v) = \Perm (v_1,v_2,\cdots, v_{d+1}) := \conv\left( \left(v_{\sigma(1)},v_{\sigma(2)},\cdots, v_{\sigma({d+1})}\right):\quad \sigma\in \fS_{d+1}\right).\]
In particular, if $\v = (1, 2, \dots, {d+1}),$ we obtain the well-known regular permutohedron. 
Note that as long as there are two different entries in $\v$, we have $\dim (\Perm(\v)) = d$. A \textbf{generic permutohedron} is any polytope of the form $\Perm(\v)$ where all the entries of $\v$ are \emph{distinct}.
\end{defn}

Recall that we have defined the braid fan $\Sigma_d$ in the introduction. Here we give a more combinatorial description of it in Lemma \ref{lem:1-1} below.
Let $\e_1,\cdots,\e_{d+1}$ be the standard basis of $\R^{d+1}$ and for each $T\in \hB_{d+1}$ we define $\e_T:=\sum_{i\in T}\e_i$ as an element in $W_d$. For any $k$-chain $T_\bullet$ of $\hB_{d+1},$ we define the corresponding \textbf{braid cone}
\[\sigma_{T_\bullet}:=\cone(\e_T: T\in T_\bullet),\]
which is $k$-dimensional. The following is well known (for a proof see \cite[Proposition 3.5]{nested}).
\begin{lem}\label{lem:1-1}
The map $T_\bullet \mapsto \sigma_{T_\bullet}$ gives a one-to-one correspondence between chains in $\CC_{d+1}$ and cones in the braid fan $\Sigma_d.$ Moreover, $k$-chains in $\CC_{d+1}$ are in bijection with $k$-dimensional cones in $\Sigma_d.$

\end{lem}

\begin{lem}\label{lem:generic}
The normal fan of any generic permutohedron is the braid fan $\Sigma_d$. 
\end{lem}


\subsection{Permutohedral variety}
For toric varieties we follow the notation and terminology of \cite{fulton}. The \textbf{permutohedral variety} $X_d$ is the toric variety associated to $\Sigma_d$ over an algebraically closed field of characteristic zero $\k$. For each $T_\bullet\in\CC_{d+1}$, its corresponding braid cone $\sigma_{T_{\bullet}}$ is associated with a subvariety $V(\sigma_{T_{\bullet}})$. These subvarieties are the \textbf{torus invariant cycles}.

For any $d\in \NN$ we define the following ring \[R_d:=\mathbb{Q}[ T\in \hB_{d+1}].\] 
For any element $i\in[d+1]$ we define the linear form $\ell_i:=\sum_{T \ni i} x_T$. 
\begin{defn}
	
	The \textbf{Chow ring of the permutohedral variety $X_d$} can be presented as
\begin{equation}\label{eq:Chowdef}
A_d \cong R_d/(I_1+I_2)
\end{equation}
where
\[
I_1 = \langle x_Tx_{T'}: \text{ for $T,T'$ incomparable}\rangle,\qquad I_2 = \langle \ell_a-\ell_b:\text{ for all $a,b\in [d+1]$} \rangle.
\]
\end{defn}

We are interested in computing the Todd class of $X_d$ in $A_d$. The following definition follows \cite[Section 5]{fulton}.
\begin{defn}
	The \textbf{Todd class of $X_d$} is the element of $A_d$ defined as
\begin{equation}\label{eq:tddef}
\td (X_d) := \prod_{T\in\hB_{d+1}} \left( \dfrac{x_T}{1-e^{-x_T}}\right),
\end{equation}
which is an element of $A_d$ by expanding each parenthesis on the right hand side as
\begin{equation}\label{eq:tdexpansion}
\dfrac{x_T}{1-e^{-x_T}} = 1+\dfrac{x_T}{2}+\sum_{i=1}^\infty \dfrac{(-1)^{i-1}B_i}{(2i)!}x_T^{2i}=1+\dfrac{x_T}{2}+\dfrac{x_T^2}{12}-\dfrac{x_T^4}{720}+\dfrac{x_T^6}{30240}+\cdots.
\end{equation}
Here $B_i$ is the $i$-th Bernoulli number. A basic fact about Chow rings is that monomials of $R_d$ of degree greater than $d$ are zero in $A_d$. Hence, the sum in \eqref{eq:tdexpansion} is finite, and thus \eqref{eq:tddef} is well-defined in $A_d$. In order to be self-contained, we will state the basic fact used above in Corollary \ref{cor:highdegree} and give an elementary proof using our computations in Section \ref{sec:gen}. 
\end{defn}

For each $T_\bullet \in \CC_{d+1},$ the class of the subvariety $V(\sigma_{T_{\bullet}})$ in $A_d$ is denoted $[V(\sigma_{T_{\bullet}})]$, and it can be represented as a square-free element in $A_d:$ 
\begin{equation}\label{eq:sqfree}
[V(\sigma_{T_{\bullet}})] = x_{T_{\bullet}}:= \prod_{T\in T_{\bullet}} x_T.
\end{equation}
We are interested in expressions for $\td(X_d)$ in terms of classes of the torus invariant cycles. In other words, we are looking for $r(T_\bullet)\in\Q$ such that
\begin{equation}\label{eq:tdexpression}
	\td(X_d) = \sum_{T_{\bullet}\in\CC_{d+1}} r_d(T_\bullet) \ x_{T_\bullet} = \sum_{T_{\bullet}\in\CC_{d+1}} r_d(T_\bullet) \ [V(\sigma_{T_{\bullet}})].
\end{equation}
We call such an expression a \textbf{square-free expression for the Todd class $\td(X_d)$ of $X_d$}.
\begin{rem}\label{rem:sqfree}
By Equation \eqref{eq:sqfree} an expression of the form \eqref{eq:tdexpression} can be obtained by finding a \emph{square-free} representation in $A_d$.
\end{rem}
Our interest in such an expression lies in the following theorem originally attributed to Danilov which is already mentioned in \S \ref{subsec:toddintro}. Here we only state it in the particular case of braid fans.


\begin{thm}[Section 5 in \cite{fulton}]\label{thm:dani}
	Let $P$ be a $d$-dimensional lattice generalized permutohedron with normal fan $\Sigma_d$. Suppose $r_d$ is a function defined on $\CC_{d+1}$ such that Equation \eqref{eq:tdexpression} holds. Using the one-to-one correspondence between chains in $\CC_{d+1}$ and cones in $\Sigma_d$ described in Lemma \ref{lem:1-1}, we can consider $r_d$ to be a function on braid cones by letting 
	\[ r_d(\sigma_{T_\bullet}) := r_d(T_\bullet).\]
	Then we have that
\begin{equation}\label{eq:mcmullen}
\Lat(P) = \sum_{F\subset P} r_{d}(\ncone(F,P)) \ \nvol(F).
\end{equation}
\end{thm}
Therefore, an equation of the form \eqref{eq:tdexpression} gives a solution to \eqref{eq:exterior} for lattice generalized permutohedra by setting $\alpha(F,P)=r_d(\ncone(F,P))$.

We are focusing on the particular case of braid fans instead of on all possible fans at the same time, so \emph{a priori} we are not looking for a Danilov function. However, we are going to require one more special property for our expressions of the form \eqref{eq:tdexpression}.

\begin{defn}\label{def:sym}
The symmetric group $\fS_{d+1}$ acts on elements of $\hB_{d+1}$ hence on the generators of the ring $R_d$. Notice that this action fixes both ideals $I_1$ and $I_2$ so that $\fS_{d+1}$ acts naturally on $A_d$ too. We say an element $f\in A_d$ is \textbf{symmetric} if $\pi\cdot f=f$ for all $\pi\in \fS_{d+1}$. 

For any $f \in A_d,$ we define its \textbf{symmetrization} to be 
\begin{equation}\label{eq:sym}
f^{\sharp}:=\dfrac{1}{(d+1)!}\sum_{\pi\in \fS_{d+1}} \pi\cdot f.
\end{equation}
(It is easy to see that $f^{\sharp}$ is symmetric.)
\end{defn}

\begin{rem}\label{rem:sym}
	In the ring $A_d$, any square-free element is of the form $\sum_{T_\bullet\in\CC_{d+1}} r(T_\bullet) \ x_{T_\bullet}$ where $r(T_\bullet) \in \mathbb{Q}$, and it is symmetric if the $r(T_\bullet)$ depends only on the \emph{size vector} of $T_\bullet$, i.e., the sequence of integers $|T_1|,|T_2|,\ldots,|T_k|$ for each $T_\bullet=(T_1,\cdots,T_k)\in\CC_{d+1}$.

\end{rem}

Recall that the BV-function $\bv$ is both a McMullen function and a Danilov function. In the case of the braid fan, we abuse notation again and consider $\bv$ a function on $\CC_{d+1}$ by letting
\[ \bv(T_\bullet) := \bv(\sigma_{T_\bullet}), \quad \forall T_\bullet \in \CC_{d+1}.\]
Then using results from \cite{bvalpha} we prove the following. 

\begin{thm}\label{thm:uniqueness}[Theorem 5.5 in \cite{bvalpha}] 
	There is a unique symmetric square-free expression for $\td(X_d).$ 
	It is given by the Berline-Vergne function: 
\begin{equation}\label{eq:bvexp}
\td(X_d) = \sum_{T_{\bullet}\in\CC_{d+1}} \bv(T_{\bullet}) [V(\sigma_{T_{\bullet}})].
\end{equation}
\end{thm}

We call the right hand side of Equation \eqref{eq:bvexp} the \textbf{Berline-Vergne expression} for the Todd class of $X_d$.

\begin{proof}
By Theorem \ref{thm:dani}, any expression for Equation \eqref{eq:tdexpression} yields an expression in the form of \eqref{eq:mcmullen} for any lattice polytope with the braid fan being its normal fan. In particular, by Lemma \ref{lem:generic}, this applies to all lattice generic permutohedra. Also, one checks that an expression is symmetric in the sense of Definition \ref{def:sym} if and only if it is symmetric in the sense of \cite[Definition 3.13]{bvalpha}. 
Hence, the conditions of \cite[Setup 4.1]{bvalpha} are met, and thus by \cite[Theorem 5.5]{bvalpha}, the conclusion of the theorem follows. 
\end{proof}

Combining the theorem with the symmetrization described in Equation \eqref{eq:sym} we get the following

\begin{prop}\label{prop:symm}
Let $f$ be any square-free expression for $\td(X_d)$ (as in Equation \eqref{eq:tdexpression}), then its symmetrization $f^{\sharp}$ is the Berline-Vergne expression for $\td(X_d)$.
\end{prop}

\section{Combinatorial Tools: Spider diagrams}\label{sec:spider}

%
%
%

In this section we develop the necessary combinatorial language that will be used to express our main formulas in Section \ref{sec:gen}.

\begin{defn}
	Let $T_\bullet\in\CC_{d+1}$ and $S\in T_\bullet$ (so $S$ is a subset of $[d+1]$). 
	
	A \textbf{spider} $\Sp = \Sp(T_\bullet, S)$ on $T_\bullet$ with head $S$ is a graph on the vertex set $T_\bullet$ with edge set $\{S,T\}$ for every $T\in T_\bullet\backslash\{S\}$.
We call $S$ the \textbf{head} and every non-head vertex a \textbf{leg}. Legs are partitioned into two subsets $L$ and $R$. The set $L$ consists of the \textbf{left legs}, the elements $T\in\T_\bullet$ such that $T\subset S$, and the set $R$ consists of the \textbf{right legs}, the elements $T\in\T_\bullet$ such that $T\supset S$.

The \textbf{size} of a spider is $|\Sp(T_\bullet,S)|:=|T_\bullet|$, the size of its vertex set. A spider of size one is called a \textbf{trivial} spider. It has no legs. (Note that number of edges in a spider $\Sp$ is $|\Sp|-1.$)
\end{defn}
\begin{ex}
We draw spiders by aligning the vertices and circling the head. In this way the sets $L$ and $R$ are visually on the left and right respectively. To save space we avoid commas, for instance $\{12\}:=\{1,2\}$.
See Figure \ref{ex:1st} for an example where $$T_\bullet =\{2\}\subset \{12\}\subset \{123\}\subset \{123456\}\subset \{12345678\}$$ and $S=\{123\}$.
\begin{figure}[ht]
\begin{tikzpicture}%
\begin{scope}

  \node at (-2,-0.5) {$\{2\}$};  
  \node at (0,-0.5) {$\{12\}$};
  \node at (2,-0.5) {$\{123\}$};
  \node at (4,-0.5) {$\{123456\}$};
  \node at (6,-0.5) {$\{12345678\}$};

  \tikzstyle{every node}=[circle, fill, draw, inner sep=1pt,scale=0.5]

	\node (i2) at (-2,0) {2};
    \node (i3) at (0,0) {3};
    \node (i4) at (2,0) {4};
    \node (i5) at (4,0) {5};
    \node (i6) at (6,0) {6};

    \draw [thick] (i2) to[out=75, in=105] (i4);
    \draw [thick] (i3) to[out=60, in=120] (i4);
    \draw [thick] (i4) to[out=60, in=120] (i5);
    \draw [thick] (i4) to[out=75, in=105] (i6);

	\draw [thick] (2,0) circle [radius=0.2];

\end{scope}

 \end{tikzpicture}
\caption{A drawn spider.}\label{ex:1st}
\end{figure}
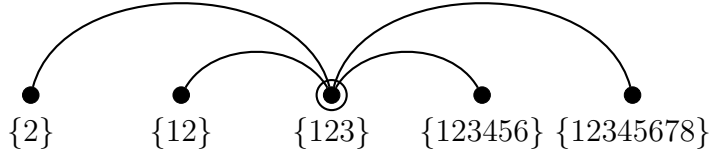
\end{ex}

Given a spider $\Sp=\Sp(T_\bullet, S),$ an \textbf{edge labeling} $\omega$ of $\Sp$ is bijection from the set of edges of $\Sp$ to the $[|\Sp|-1] =\{1, 2, \dots, |\Sp|-1\}.$ We say an edge labeling $\omega$ of $\Sp$ is \textbf{natural} if $\omega(\{S,T\})>\omega(\{S,T'\})$ whenever $S\subset T\subset T'$ or $T'\subset T\subset S$. We use notation $(\Sp, \omega)$ to indicate a spider $\Sp$ with a natural edge labeling $\omega$. 
\begin{ex}
We draw a natural edge labeling on the spider of Figure \ref{ex:1st}.
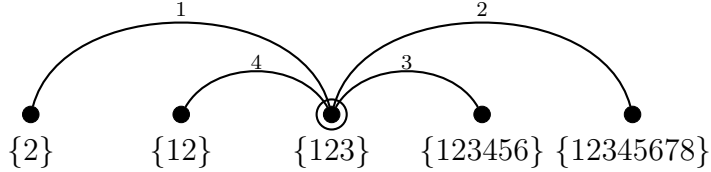
\begin{figure}[ht]
\begin{tikzpicture}%
\begin{scope}

  \node at (-2,-0.5) {$\{2\}$};  
  \node at (0,-0.5) {$\{12\}$};
  \node at (2,-0.5) {$\{123\}$};
  \node at (4,-0.5) {$\{123456\}$};
  \node at (6,-0.5) {$\{12345678\}$};

  \node at (0,1.4) {\tiny{1}};
  \node at (1,.7) {\tiny{4}};
  \node at (3,.7) {\tiny{3}};
  \node at (4,1.4) {\tiny{2}};

  \tikzstyle{every node}=[circle, fill, draw, inner sep=1pt,scale=0.5]

	\node (i2) at (-2,0) {2};
    \node (i3) at (0,0) {3};
    \node (i4) at (2,0) {4};
    \node (i5) at (4,0) {5};
    \node (i6) at (6,0) {6};

    \draw [thick] (i2) to[out=75, in=105] (i4);
    \draw [thick] (i3) to[out=60, in=120] (i4);
    \draw [thick] (i4) to[out=60, in=120] (i5);
    \draw [thick] (i4) to[out=75, in=105] (i6);

	\draw [thick] (2,0) circle [radius=0.2];

\end{scope}

 \end{tikzpicture}
\caption{A drawn spider with a natural edge labeling.}
\end{figure}
\end{ex}

\begin{notn}\label{notn:label}
A left leg will be labeled as $T_i^L$ if it is the $i$-th \emph{smallest} vertex among all left legs, and a right leg will be labeled by $T_j^R$ if it is the $j$-th \emph{largest} vertex among all right legs.
If there are no left legs, we give the head vertex $S$ an additional label $T^L_1$ similarly, if there are no right legs, we give the head vertex $S$ an additional label $T^R_1$. 
\end{notn}

\begin{ex}

Consider the spider with chain $T_\bullet = \{12\}\subset \{123\}\subset \{123456\}\subset \{12345678\}$ and head $S=\{12\}$. In Figure \ref{fig:singlespider} we have labeled the spider according to Notation \ref{notn:label}. Note that the head $S$ also receives the label $T^L_1$.
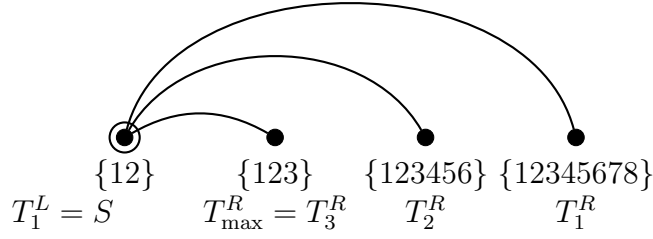
\begin{figure}[ht]
\begin{tikzpicture}%
\begin{scope}

  \node[left] at (0,-1) {$T^L_1=S$};
  \node at (2,-1) {$T^R_{\max}=T^R_{3}$};
  \node at (4,-1) {$T^R_{2}$};
  \node at (6,-1) {$T^R_{1}$};
  
  \node at (0,-0.5) {$\{12\}$};
  \node at (2,-0.5) {$\{123\}$};
  \node at (4,-0.5) {$\{123456\}$};
  \node at (6,-0.5) {$\{12345678\}$};

  \tikzstyle{every node}=[circle, fill, draw, inner sep=1pt,scale=0.5]


    \node (i3) at (0,0) {3};
    \node (i4) at (2,0) {4};
    \node (i5) at (4,0) {5};
    \node (i6) at (6,0) {6};

    \draw [thick] (i3) to[out=30, in=150] (i4);
    \draw [thick] (i3) to[out=60, in=120] (i5);
    \draw [thick] (i3) to[out=75, in=105] (i6);

	\draw [thick] (0,0) circle [radius=0.2];

\end{scope}

 \end{tikzpicture}
 \caption{A spider with the head having two different labels.}
 \label{fig:singlespider}
\end{figure}
\end{ex}

\begin{defn}\label{def:sd}
Let $T_\bullet\in\CC_{d+1}$ be a chain. A \textbf{spider diagram} $\D$ on $T_\bullet$ consist of a partition of $T_\bullet$ into $k$-disjoint intervals $T_{1,\bullet},\cdots,T_{k,\bullet}$ together with a spider $\Sp_i:=\Sp(T_{i,\bullet},S_i)$ on each interval. The set of heads $S_i$ is called the \textbf{head set} of $\D$. Note that the head set is always a chain $S_\bullet\in\CC_{d+1}$. For each $i$, let $m_i=|\Sp_i|$, and $L_i$ and $R_i$ be set of left and right legs respectively. The vector $\m:=(m_1,\cdots,m_k)$ is the \textbf{length vector} of $\D$. The \textbf{size} of $\D$ is $|\D|:=\sum m_i$, the size of its vertex set. An \textbf{ordered spider diagram} $\D$ is a spider diagram in which additionally we have a natural edge labeling $\omega_i$ on each spider $\Sp_i$. 

A pair $(S_\bullet,\m)$ is \textbf{admissible} if $S_\bullet\in\CC^k_{d+1}$, $\m\in\ZZ^k$. Such a pair is \textbf{d-admissible} if furthermore $\sum_{i=1}^k{m_i}\leq d$. 
Let $\spd(S_\bullet,\m)$ (respectively $\lspd(S_\bullet,\m)$) be the set of spider diagrams (respectively ordered spider diagrams) with head set $S_\bullet$ and length vector $\m$.
\end{defn}

\begin{rem}
Notice that if $(S_\bullet,\m)$ is not d-admissible, that is, $\sum_{i=1}^k{m_i}> d$, then $\spd(S_\bullet,\m)$ and $\lspd(S_\bullet,\m)$ are empty.
\end{rem}

\begin{notn}\label{notn:triple}
In a spider diagram $\D$ the legs are now triply indexed: the element $T^P_{i,j}$ with $P\in\{L,R\}$ is the $j$th smallest/largest on the side $L/R$ of the $i$-th spider $\Sp_i$.

We also let $l_i$ be the number of left legs of $\Sp_i$ and $r_i$ be the number of right legs of $\Sp_i.$ Hence, $T^L_{i, l_i}$ is the largest vertex among all left legs and $T^R_{i, r_i}$ is the largest vertex among all right legs.

\end{notn}
\begin{ex}
In Figure \ref{fig:spiderex} we show an ordered spider diagram of two spiders where vertices are labeled according to Notation \ref{notn:triple}. This diagram has size 10 and length vector $(6,4)$.  We have also included a natural edge labeling $\omega.$
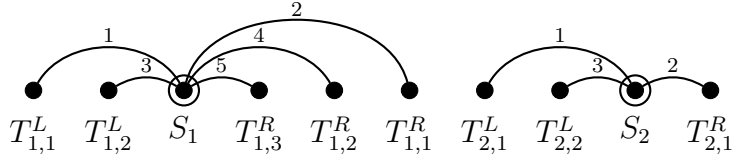
\begin{figure}[ht]\label{ex:spiderdiagram}

\begin{tikzpicture}%
\begin{scope}
	\node [below] at (1,-0.2) {$T^L_{1,1}$};
	\node [below] at (2,-0.2) {$T^L_{1,2}$};
	\node [below] at (3,-0.2) {$S_1$};
	\node [below] at (4,-0.2) {$T^R_{1,3}$};
	\node [below] at (5,-0.2) {$T^R_{1,2}$};
	\node [below] at (6,-0.2) {$T^R_{1,1}$};
	\node [below] at (7,-0.2) {$T^L_{2,1}$};
	\node [below] at (8,-0.2) {$T^L_{2,2}$};
	\node [below] at (9,-0.2) {$S_2$};
	\node [below] at (10,-0.2) {$T^R_{2,1}$};
	
	\node [above] at (2,0.5) {\tiny{1}};
	\node [above] at (2.5,0.1) {\tiny{3}};
	\node [above] at (3.5,0.1) {\tiny{5}};
	\node [above] at (4,0.5) {\tiny{4}};
	\node [above] at (4.5,0.85) {\tiny{2}};
	
	\node [above] at (8,0.5) {\tiny{1}};
	\node [above] at (8.5,0.1) {\tiny{3}};
	\node [above] at (9.5,0.1) {\tiny{2}};
	
  \tikzstyle{every node}=[circle, fill, draw, inner sep=1pt,scale=0.5]

    \node (i1) at (1,0) {1};
    \node (i2) at (2,0) {2};
    \node (i3) at (3,0) {3};
    \node (i4) at (4,0) {4};
    \node (i5) at (5,0) {5};
    \node (i6) at (6,0) {6};
    \node (i7) at (7,0) {7};
    \node (i8) at (8,0) {8};
    \node (i9) at (9,0) {9};
    \node (j1) at (10,0) {1};

    \draw [thick] (i3) to[out=150, in=30] (i2);
    \draw [thick] (i3) to[out=120, in=60] (i1);
    \draw [thick] (i3) to[out=30, in=150] (i4);
    \draw [thick] (i3) to[out=60, in=120] (i5);
    \draw [thick] (i3) to[out=75, in=105] (i6);
    \draw [thick] (i9) to[out=150, in=30] (i8);
    \draw [thick] (i9) to[out=120, in=60] (i7);
    \draw [thick] (i9) to[out=30, in=150] (j1);   
    
    \draw [thick] (3,0) circle [radius=0.2]; 
	\draw [thick] (9,0) circle [radius=0.2];


\end{scope}

 \end{tikzpicture}
 \caption{A spider diagram of two spiders together with a natural edge labeling.}
 
 \label{fig:spiderex}
\end{figure}
\end{ex}

For our formulas in the next section we need to define the weight of a spider diagram.
\begin{defn}\label{def:weight}
Let $T_\bullet\in\CC_{d+1}$ be a chain and $\D$ a spider diagram on $T_\bullet$ with $k$ spiders. 
We define the \textbf{internal weight} of a single spider $\Sp_i$ as
\begin{equation}
	\intwt(\Sp_i):=\left(\prod_{j>1}\dfrac{|T^L_{i,j}-T^L_{i,j-1}|}{|S_i-T^L_{i,j-1}|}\right)\left(\prod_{j>1}\dfrac{|T^R_{i,j-1}-T^R_{i,j}|}{|T^R_{i,j-1}-S_i|}\right), 
\end{equation} 
(note that the internal weights of a trivial spider is 1) and the \textbf{boundary weight} of the diagram $\D$ as
\begin{equation}
\bdwt(\D):=\dfrac{|T^L_{1,1}-\emptyset|}{|S_1-\emptyset|}\left(\prod_{i=2}^{k}\dfrac{|T^L_{i,1}-T^R_{i-1,1}|}{|S_i-S_{i-1}|}\right)\dfrac{|[d+1]-T^R_{k,1}|}{|[d+1]-S_{k}|}.
\end{equation}

The \textbf{weight} of a spider diagram $\D$ is defined as
\[
\wt(\D):=\bdwt(\D)\prod_{i=1}^{k} \intwt(\Sp_i).
\]
\end{defn}
\begin{ex}
The weight of the spider diagram $\D$ depicted in Figure \ref{fig:spiderex} is
\begin{align*}
\underbrace{\left(\dfrac{|T_{1,1}^L-\emptyset|}{|S_1-\emptyset|}\cdot\dfrac{|T_{2,1}^L-T_{1,1}^R|}{|S_2-S_{1}|}\cdot\dfrac{|[d+1]-T_{2,1}^R|}{|[d+1]-S_2|}\right)}_{\bdwt(\D)}
&\underbrace{\left(\dfrac{|T_{1,2}^L-T_{1,1}^L|}{|S_1-T_{1,1}^L|}\right)\left(\dfrac{|T_{1,1}^R-T_{1,2}^R|}{|T_{1,1}^R-S_1|}\cdot\dfrac{|T_{1,2}^R-T_{1,3}^R|}{|T_{1,2}^R-S_1|}\right)}_{\intwt(\Sp_1)}\\
&\times \underbrace{\left(\dfrac{|T_{2,2}^L-T_{2,1}^L|}{|S_2-T_{2,1}^L|}\right)\cdot 1}_{\intwt(\Sp_2)}.
\end{align*}

\end{ex}

%


\section{Computations in the Chow ring $A_d$}\label{sec:gen}

The main goal of this section is to express any element in $A_d$ as a sum of square-free monomials. We start by treating squares.

\begin{lem}\label{lem:nonsym}
Let $S\in\hB_{d+1}$. Choose $a\in S$, $b\notin S$ then
\begin{equation}\label{eq:nonsym}
x_S^2 = -\sum_{\substack{T\subset S \\ a\in T}}x_Tx_S -\sum_{\substack{S\subset T \\ b\notin T}} x_Sx_T.
\end{equation}
\end{lem}
\begin{proof}
Using the relation $\l_a-\l_b\in I_2$, we get $x_S(\l_a-\l_b)=0$ in $A_d$. Hence, 
\begin{equation}
x_S\left(\sum_{a\in T}x_T-\sum_{b\in T}x_T\right) = 0.
\end{equation}
The relations in $I_1$ imply that $x_S x_T=0$ for any $T$ that is neither $T\subseteq S$ nor $S\subseteq T$. Thus we expand the above equation to get.

\begin{equation}
x_S^2+\sum_{ \substack{T\subset S\\ a\in T}} x_Tx_S + \sum_{\substack{S\subset T\\ a\in T}} x_Sx_T - \sum_{ \substack{T\subset S\\ b\in T}} x_Tx_S - \sum_{\substack{S\subset T\\ b\in T}} x_Sx_T = 0.
\end{equation}
The fourth term is zero since the condition is vacuous. In the third term notice that the condition $a\in T$ is redundant. After canceling terms from the second and fourth sums, everything reduces to
\[
x_S^2+\sum_{\substack{T\subset S \\ a\in T}}x_Tx_S + \sum_{\substack{S\subset T \\ b\notin T}}x_Sx_T = 0.
\]
By solving for $x_S^2$ we get Equation \eqref{eq:nonsym}.
\end{proof}

The square-free expression obtained in Lemma \ref{lem:nonsym} is not symmetric. To adjust this, we average over all possibilities.
\begin{lem}\label{lem:squares}
Let $S_1,\cdots,S_k$ be a k-chain in $\hB_{d+1}$ and $S_l$ a set in the chain. We have the following equality in $A_d$:
\begin{align}
x_{S_1}\cdots x_{S_{l-1}}x_{S_l}^2x_{S_{l+1}}\cdots x_{S_k} =& \ -\sum_{S_{l-1}\subset T\subset S_l} \dfrac{|T-S_{l-1}|}{|S_l-S_{l-1}|} x_{S_1}\cdots x_{S_{l-1}}x_Tx_{S_l}x_{S_{l+1}}\cdots x_{S_k} \label{eq:squares}\\
& \ -\sum_{S_{l}\subset T\subset S_{l+1}} \dfrac{|S_{l+1}-T|}{|S_{l+1}-S_{l}|} x_{S_1}\cdots x_{S_{l}}x_Tx_{S_{l+1}}x_{S_{l+2}}\cdots x_{S_k}\nonumber
\end{align}
By convention, we let $S_0=\emptyset$ and $S_{k+1}=[d+1]$.
\end{lem}
\begin{proof}
We expand $x_{S_l}^2$ as in Lemma \ref{lem:nonsym} using all possible pairs $(a,b)\in (S_{l}-S_{l-1})\times (S_{l+1}-S_{l})$ and take the average. We obtain
\begin{align*}
& \ x_{S_1}\cdots x_{S_{l-1}}x_{S_l}^2x_{S_{l+1}}\cdots x_{S_k} \\
=& \ -\frac{1}{|S_l-S_{l-1}|\cdot |S_{l+1}-S_l|} \left( \sum_{\substack{a \in S_l- S_{l-1} \\ b \in S_{l+1}-S_l}} \sum_{\substack{S_{l}\subset T\subset S_{l+1} \\ a\in T}} x_{S_1}\cdots x_{S_{l-1}}x_Tx_{S_l}x_{S_{l+1}}\cdots x_{S_k}\right) \\
& \ - \frac{1}{|S_l-S_{l-1}|\cdot |S_{l+1}-S_l|} \left(\sum_{\substack{a \in S_l- S_{l-1} \\ b \in S_{l+1}-S_l}} \sum_{\substack{S_l\subset T\subset S_{l+1} \\ b\notin T}} x_{S_1}\cdots x_{S_{l}}x_Tx_{S_{l+1}}x_{S_{l+2}}\cdots x_{S_k} \right)
\end{align*}
For the first term in the above expression, one sees that a set $T$ contributes to the summation if and only if $S_{l-1}\subset T\subset S_l$, and for any $T: S_{l-1}\subset T\subset S_l,$ it appears if and only if $a \in T-S_{l-1}$ which can be paired with any $b \in S_{l+1}-S_l,$ and thus it appears exactly $|T-S_{l-1}|\cdot |S_{l+1}-S_l|$ times. Hence, the first term above agrees with the first term on the right hand side of \eqref{eq:squares}. Similarly, we can show that the second term above coincides with the second term on the right hand side of \eqref{eq:squares}.
\end{proof}
\begin{rem}\label{rem:zerocase}
	The two sums on the right hand side of Equation \eqref{eq:squares} can be over empty index sets simultaneously, in which case the monomial on the left hand side of \eqref{eq:squares} is equal to zero. In particular, notice that this happens when $(S_{l}-S_{l-1})\times (S_{l+1}-S_{l})$ is a singleton.
\end{rem}


Repeated use of this lemma allows us to expand any monomial in $A_d$ as a sum of squarefree monomials. First some more notation.

\begin{defn}
For a spider diagram $\D$ we define $\sgn(\D):=(-1)^{|\D|-k}$. The number $|\D|-k$ is equal to the total number of legs. Also we let $x_\D:=x_{T_\bullet}$, where $T_\bullet$ is the vertex set of $\D$.
\end{defn}

\begin{prop}\label{prop:expanding}
Let $S_\bullet\in\CC_{d+1}^k$ be a $k$-chain and $\m = (m_1,\cdots,m_k)$ a vector of positive integers. 
Recall from Definition \ref{def:sd} that $\lspd(S_\bullet,\m)$ is the set of ordered spider diagrams where each spider $\Sp_i$ has head $S_i$, size $m_i$, and an edge labeling $\omega_i$.
Then we have the following equality in the ring $A_d$:
\begin{equation}\label{eq:expansion}
x_{S_\bullet}^\m:=\prod_{i=1}^k x_{S_i}^{m_i} = \sum_{\D\in\lspd(S_\bullet,\m)}\sgn(\D)\wt(\D)x_\D.
\end{equation}
\end{prop}

\begin{proof}
We prove by induction on $|\m|=\sum m_i$. 
	The base case is $|\m| = k$ which happens exactly when $m_i=1$ for all $i$. In this case we have the square-free monomial $\prod_{i=1}^k x_{S_i}$ and $\lspd(S_\bullet,\m)$ consists of a single diagram $\D_0$ with $k$ trivial spiders. One sees that $\sgn(\D_0)=1$ and $\wt(\D_0)=1$ so Equation \eqref{eq:expansion} is trivially true.

	We proceed to the induction step. Suppose $N > k$ is a positive integer, and for any $\m=(m_1, \dots, m_k)$ with $|\m| = \sum m_i < N,$ we have the equality \eqref{eq:expansion} holds in $A_d.$
Now we assume $\m=(m_1, \dots, m_k)$ is a vector of positive integers satisfying $|\m| = \sum m_i =N.$
Let $j=\max\{i:m_i>1\}$ and $\m':=(m_1,\cdots,m_j-1,\cdots)$. By the induction hypothesis, we have the following equality in $A_d:$
\begin{equation}\label{eq:induction}
x_{S_j}^{-1}\prod_{i=1}^k x_{S_i}^{m_i}=\sum_{\D'\in\lspd(S_\bullet,\m')}\sgn(\D')\wt(\D')x_{\D'}.
\end{equation}
We see that in order to show that \eqref{eq:expansion} holds for $\m$, it is enough to show
\begin{equation}\label{eq:needtoshow}
\sum_{\D'\in\lspd(S_\bullet,\m')}\sgn(\D')\wt(\D')x_{S_j}x_{\D'} = \sum_{\D\in\lspd(S_\bullet,\m)}\sgn(\D)\wt(\D)x_\D.
\end{equation}


For the rest of the proof, we will use notations developed in Section 3, in particular recall those given in Notation \ref{notn:triple}.
In order to prove \eqref{eq:needtoshow}, we construct a \textbf{pruning} map $\prune$ from $\lspd(S_\bullet,\m)$ to $\lspd(S_\bullet, \m')$ in the following way: Let $\D \in \lspd(S_\bullet,\m)$, where $\Sp_j$ is the $j$th spider in $\D$ together with a natural edge labeling $\omega_j$. Suppose $T$ is the leg in $\Sp_j$ such that $\omega_j\left(\{S_j,T\}\right)$ has the largest edge label in $\Sp_j$. (Note that $T$ is either the closest leg/vertex $T^L_{j,l_j}$ on the left of $S_j$ or the closest leg/vertex $T^R_{j,r_j}$ on the right of $S_j$.) Then we define $\prune(\D)$ to be the ordered spider diagram obtained from $\D$ by removing $T$. For example, if $\D$ is the ordered spider diagram in Figure \ref{fig:spiderex}, then $\prune(\D)$ is the one depicted in Figure \ref{fig:prune}.
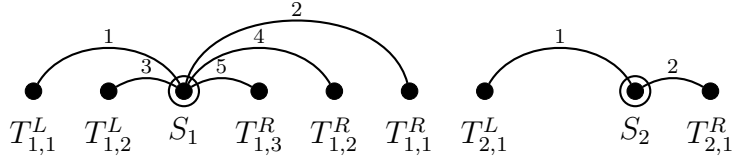
\begin{figure}[ht]
\begin{tikzpicture}%
\begin{scope}
	\node [below] at (1,-0.2) {$T^L_{1,1}$};
	\node [below] at (2,-0.2) {$T^L_{1,2}$};
	\node [below] at (3,-0.2) {$S_1$};
	\node [below] at (4,-0.2) {$T^R_{1,3}$};
	\node [below] at (5,-0.2) {$T^R_{1,2}$};
	\node [below] at (6,-0.2) {$T^R_{1,1}$};
	\node [below] at (7,-0.2) {$T^L_{2,1}$};
	\node [below] at (9,-0.2) {$S_2$};
	\node [below] at (10,-0.2) {$T^R_{2,1}$};
	
	\node [above] at (2,0.5) {\tiny{1}};
	\node [above] at (2.5,0.1) {\tiny{3}};
	\node [above] at (3.5,0.1) {\tiny{5}};
	\node [above] at (4,0.5) {\tiny{4}};
	\node [above] at (4.5,0.85) {\tiny{2}};
	
	\node [above] at (8,0.5) {\tiny{1}};
	\node [above] at (9.5,0.1) {\tiny{2}};
	
  \tikzstyle{every node}=[circle, fill, draw, inner sep=1pt,scale=0.5]

    \node (i1) at (1,0) {1};
    \node (i2) at (2,0) {2};
    \node (i3) at (3,0) {3};
    \node (i4) at (4,0) {4};
    \node (i5) at (5,0) {5};
    \node (i6) at (6,0) {6};
    \node (i7) at (7,0) {7};
    \node (i9) at (9,0) {9};
    \node (j1) at (10,0) {1};

    \draw [thick] (i3) to[out=150, in=30] (i2);
    \draw [thick] (i3) to[out=120, in=60] (i1);
    \draw [thick] (i3) to[out=30, in=150] (i4);
    \draw [thick] (i3) to[out=60, in=120] (i5);
    \draw [thick] (i3) to[out=75, in=105] (i6);
    \draw [thick] (i9) to[out=120, in=60] (i7);
    \draw [thick] (i9) to[out=30, in=150] (j1);   
    
    \draw [thick] (3,0) circle [radius=0.2]; 
	\draw [thick] (9,0) circle [radius=0.2];


\end{scope}

 \end{tikzpicture}
\caption{An example of the pruning function.}
\label{fig:prune}
\end{figure}

As sets we have $\lspd(S_\bullet,\m)=\coprod_{\D'\in \lspd(S_\bullet,\m')} \prune^{-1}(\D')$, 
so we can rewrite the right hand side of \eqref{eq:needtoshow} as
\[  \sum_{\D'\in\lspd(S_\bullet,\m')}\sum_{\D\in\prune^{-1}(\D')}\sgn(\D)\wt(\D)x_\D \quad \text{or} \quad \sum_{\D'\in\lspd(S_\bullet,\m')} \sum_{\substack{\D\in\lspd(S_\bullet,\m)\\ \prune(\D)=\D'}} \sgn(\D)\wt(\D)x_\D.\]
Hence, we can reduce the problem of proving \eqref{eq:needtoshow} to proving that for every $\D' \in \lspd(S_\bullet, \m'),$ 
\begin{equation}\label{eq:needtoshow2}
\sgn(\D')\wt(\D')x_{S_j}x_{\D'}  = \sum_{\substack{\D\in\lspd(S_\bullet,\m)\\ \prune(\D)=\D'}} \sgn(\D)\wt(\D)x_\D.
\end{equation}
We now apply Lemma \ref{lem:squares} to rewrite $x_{S_j} x_{\D'}$. One notices that the two summations appear on the right side of \eqref{eq:squares} correspond to reversing the ``pruning'' operation by adding a left or a right leg back. Hence,
\begin{equation}\label{eq:application}
x_{S_j}\cdot x_{\D'} = -\sum_{\substack{\D\in\lspd(S_\bullet,\m)\\ \prune(\D)=\D'}} c_j(\D) x_\D,
\end{equation}
where 
\[ c_j(\D) := \begin{cases}
		\dfrac{|T^L_{j,l_j}-T^L_{j,l_j-1}|}{|S_j-T^L_{j,l_j-1}|}, & \quad 
		\text{if the left leg $T^L_{j,l_i}$ is removed when pruning $\D$;}  \\
		\dfrac{|T^R_{j,r_j-1}-T^R_{j,r_j}|}{|T^R_{j,r_j-1}-S_j|}, & \quad 
		\text{if the right leg $T^R_{j,r_i}$ is removed when pruning $\D$.}  \\
	\end{cases}
\]
Also, if $l_j=1$ we let $T^L_{j, l_j-1} = T^R_{j-1, 1}$, and if $r_j=1$ we let $T^R_{j, r_j-1} := T^L_{j+1,1}.$ (Note that if $j=1,$ we consider $T^R_{j-1,1} = S_0 = \emptyset$; likewise, if $j=k,$ we consider $T^L_{k+1,1} = S_{k+1} = [d+1].$) 
Plugging \eqref{eq:application} into the left hand side of \eqref{eq:needtoshow2}, we obtain
\begin{equation}\label{eq:application2}
\sgn(\D')\wt(\D')x_{S_j}x_{\D'} = \sum_{\substack{\D\in\lspd(S_\bullet,\m)\\ \prune(\D)=\D'}} - \sgn(\D')\wt(\D')c_j(\D) x_\D.
\end{equation}
Comparing it with the right hand side of \eqref{eq:needtoshow2} and observing that $-\sgn(\D')=\sgn(\D)$ when $\D' =\prune(\D)$, one sees that the proof is completed if we can prove that for any $\D \in \lspd(S_\bullet,\m)$, if $\D'=\prune(\D)$, then
\begin{equation}\label{eq:comparison2}
	\wt(\D')c_j(\D)=\wt(\D).
\end{equation}

Suppose $T$ is the leg that is removed when we ``prune'' $\D$ to obtain $\D'.$ We will only consider the case when $T= T^L_{j,l_j}$ is a left leg of $\Sp_j.$ (The case when $T$ is a right left of $\Sp_j$ can be proved analogously.) 
It is straightforward to check from the definitions of weights that $c_j(\D')\wt(\D')=\wt(\D)$ when $l_j >1,$ i.e., $T^L_{j,l_j}$ is not the only left leg of $\Sp_j$. Indeed, in this case the boundary weights of $\D$ and $\D'$ are the same and the internal weights of $\D$ and $\D'$ differ exactly by a factor $c_j(\D)$ so \eqref{eq:comparison2} holds. 

Suppose $l_j=1.$ Thus $T^L_{j,1}$ is the only left leg of $\Sp_j$ (which is the $j$th spider in $\D$), and the $j$-th spider in $\D'$ has no left legs. Then the internal weights of $\D$ and $\D'$ are the same, whereas the boundary weights are different. Comparing $\bdwt(\D')$ and $\bdwt(\D)$, we see all but one factors in their expression are the same. The different factors $\bdwt(\D')$ and $\bdwt(\D)$ are
\[ \dfrac{|S_j-T^R_{j-1,1}|}{|S_j-S_{j-1}|} \text{ and } \frac{|T^L_{j,1}-T^R_{j-1,1}|}{|S_j-S_{j-1}|},\]
respectively.
Since
\[ c_j(\Delta) =\dfrac{|T^L_{j,l_j}-T^L_{j,l_j-1}|}{|S_j-T^L_{j,l_j-1}|} 
= \dfrac{|T^L_{j,1}-T^R_{j-1,1}|}{|S_j-T^R_{j-1,1}|}.\]
we conclude that $\bdwt(\D') \c_j(\D) =\bdwt(\D)$. Therefore, \eqref{eq:comparison2} follows, completing the proof.
\end{proof}

\begin{cor}\label{cor:highdegree}
Monomials in $R_d$ of degree larger than $d$ vanish in $A_d$.
\end{cor}
\begin{proof}
	Because of the relations in the ideal $I_1$, we only need to consider monomials of the form $x_{S_\bullet}^\m$ where $(S_\bullet,\m)$ is an admissible pair. Applying Equation \eqref{eq:expansion} we get an empty sum on the right if $\sum_i m_i>d$, hence a monomial $x_{S_\bullet}^\m$ is nonzero only if $(S_\bullet,\m)$ is d-admissible.
\end{proof}
The proof of Proposition \ref{prop:expanding} prunes one leg from a spider diagram at the time until we obtain a spider diagram consisting only of trivial spiders. The natural edge labelings are used to keep track of the order in which we remove the legs. However, one notices that if $\D_1, \D_2 \in \lspd(S_\bullet,\m)$ are on the same spider diagram $\D$ with two different natural edge labelings, 
then 
\[ \sgn(\D_1)\wt(\D_1)x_{\D_1} = \sgn(\D)\wt(\D)x_\D = \sgn(\D_2)\wt(\D_2)x_{\D_2}.\]
We have the following immediate consequence to Proposition \ref{prop:expanding}:

\begin{cor}\label{cor:expanding}
Let $S_\bullet\in\CC_{d+1}^k$ be a $k$-chain and $\m= (m_1,\cdots,m_k)$ a vector of positive integers. Recall from Definition \ref{def:sd} that $\spd(S_\bullet,\m)$ is the set of \emph{unordered} spider diagrams where each spider $\Sp_i$ has head $S_i$ and size $m_i$. Then we have the following equality in the ring $A_d$:
\begin{equation}\label{eq:expansionsimpl}
x_{S_\bullet}^\m = \sum_{\D\in\spd(S_\bullet,\m)}\Bi(\D)\sgn(\D)\wt(\D)x_\D,
\end{equation}
where
\[\Bi(\D):=\prod_{i=1}^k \binom{|\Sp_i|-1}{|L_i|,|R_i|}\]
counts the number of natural edge labelings on $\D.$
\end{cor}

\section{Formulas for the Berline-Vergne function}\label{sec:bv}
With Corollary \ref{cor:expanding} in hand, we can now write down square-free expressions of any element in $A_d$, in particular of the Todd class. First we need one more piece of notation. Any monomial $\u\in R_d$ that is non-zero in $A_d$ is of the form $x_{S_\bullet}^\m$ for some d-admissible pair $(S_\bullet,\m)$.
\begin{defn}
Expand the Todd class by plugging Equation \eqref{eq:tdexpansion} into Equation \eqref{eq:tddef} to obtain:
\begin{equation}\label{eq:masterbegin}
\td(X_d) = \sum_{(S_\bullet,\m)}\Tdcoeff(S_\bullet,\m)\cdot x_{S_\bullet}^\m,\quad \Tdcoeff(S_\bullet,\m)\in\mathbb{Q},
\end{equation}
where the summation is over all d-admissible pairs $(S_\bullet,\m)$. For $\D\in\spd(S_\bullet,\m)$ define $\Tdcoeff(\D):=\Tdcoeff(S_\bullet,\m)$.

\end{defn}
\begin{ex}
Let $d\geq 7$, $S_\bullet=(S_1,S_2,S_3)$ be any chain of length three and $\m=(2,4,1)$. The monomial $x_{S_1}^2x_{S_2}^4x_{S_3}^1$ appears in Equation \eqref{eq:tdexpansion} with coefficient
\[
\left(\frac{1}{12}\right)\left(-\frac{1}{720}\right)\left(\frac{1}{2}\right) = -\frac{1}{17280},
\]
hence $\Tdcoeff(S_\bullet,(2,4,1))=-1/17280$.
\end{ex}

\begin{rem}\label{rem:Tsymm}
Notice that by definition $\Tdcoeff$ only depends on the length and entries of $\m$.
\end{rem}

Now we can present our main result.
\begin{thm}\label{thm:main}
Let $X_d$ be the permutohedral variety. Its Todd class has the following representation in terms of toric invariant cycles:
\begin{equation}\label{eq:mainexpansion}
\td(X_d)= \sum_{T_{\bullet}\in\CC_{d+1}} \alpha(T_{\bullet}) [V(\sigma_{T_{\bullet}})],
\end{equation}
with coefficients $\alpha(T_\bullet)$ given by following explicit combinatorial formula
\begin{equation}\label{eq:masterformula}
\alpha(T_{\bullet}) = \sum_{\D\in\mathcal{D}(T_\bullet)} \Tdcoeff(\D)\Bi(\D)\sgn(\D)\wt(\D),
\end{equation}
where $\mathcal{D}(T_\bullet)$ is the set of all spider diagrams on $T_\bullet$. Furthermore, we have that $\alpha(T_\bullet)=\bv(T_\bullet)$, where $\bv(\cdot)$ is the Berline-Vergne function.
\end{thm}

\begin{proof}
We start with the expansion in Equation \eqref{eq:masterbegin}. Then we expand each $ x_{S_\bullet}^\m$ by using Corollary \ref{cor:expanding}. We obtain
\begin{equation}\label{eq:masterproof1}
\td(X_d) = \sum_{\D}\Tdcoeff(\D)\Bi(\D)\sgn(\D)\wt(\D)x_\D,
\end{equation}
Where the sum is over \emph{all} possible spider diagrams, to be more precise over $\bigcup\spd(S_\bullet,\m)$ where the union is over all d-admissible pairs.  We can rearrange the sum as follows
\begin{equation}\label{eq:masterproof2}
\td(X_d) = \sum_{T_\bullet\in\CC_{d+1}}x_{T_\bullet}\left(\sum_{\D\in\mathcal{D}(T_\bullet)}\Tdcoeff(\D)\Bi(\D)\sgn(\D)\wt(\D)\right).
\end{equation}
We proceed to show that the expression obtained in \eqref{eq:masterproof2} is symmetric. By Remark \ref{rem:sym} it is enough to prove that two different chains with the same size vector have the same coefficient.

Consider two chains $T_\bullet$ and $T'_\bullet$ with same size vector and fix a bijection $\phi$ on $[d+1]$ that simultaneously bijects $T_i$ with $T'_i$ for all relevant $i$. The function $\phi$ also induces a natural bijection (which abusing notation we call $\phi$ also) $\phi:\mathcal{D}(T_\bullet)\to\mathcal{D}(T'_\bullet)$. By Remark \ref{rem:Tsymm} $\Tdcoeff(\D)=\Tdcoeff(\phi(\D))$. Also $\sgn(\D)=\sgn(\phi(\D))$ since they both have the same number of legs. Finally $(\Bi(\D),\wt(\D))=(\Bi(\phi(\D)),\wt(\phi(\D)))$ since by definition $\Bi$ and $\wt$ depend only on the sizes involved.

We have thus proved that \eqref{eq:masterproof2} is a symmetric square-free expression of $\td(X_d)$. By Theorem \ref{thm:uniqueness} there is only one such expression, the one given by the Berline-Vergne function, hence we obtain the last part of the theorem.
\end{proof}

To end this section we count the number of terms appearing in Equation \eqref{eq:masterformula}.

\begin{prop}
The number of terms in Equation \eqref{eq:masterformula} is exponential.
\end{prop}

\begin{proof}
We are looking for $h(n)$ be the number of spider diagrams $\D$ on $[n]$ such that $|\Sp|\in\{1,2,4,6,\cdots\}$ for each $\Sp\in \D$. Its generating function is
\begin{equation}\label{eq:genfct1}
\sum_{n=z1}^{\infty} h(n)z^n = \sum_{k=1}^{\infty}\left(x+2x^2+4x^4+6x^6+\cdots\right)^k.
\end{equation}
Observe that $x+2x^2+4x^4+6x^6+\cdots = x + \frac{1}{2}\left(\frac{x}{(1-x)^2}+\frac{(-x)}{(1-(-x)^2)}\right)=\frac{x(x^4-2x^2+2x+1)}{(x^2-1)^2}$. Plugging into each term in the right of \eqref{eq:genfct1} and using the geometric series formula we obtain
\begin{equation}\label{eq:genfct2}
\sum_{n=z1}^{\infty} h(n)z^n = -\dfrac{z(z^4-2z^2+2z+1)}{z^5-z^4-2z^3+4z^2+z-1}=z+3z^2+5z^3+15z^4+29z^5+\cdots.
\end{equation}
The conclusion follows since the coefficients of a rational function are asymptotically the powers of the largest root of the denominator \cite[Theorem 4.1.1]{ec2}, which in this case is $\approx 1.602$.
\end{proof}

\section{Examples in low (co)dimension}\label{sec:ex}
In this section we explicitly compute Equation \eqref{eq:masterformula} for chains of small length.

\begin{prop}[Codimension 2 cones.]
Let $(T_1,T_2)\in\CC^2_{d+1}$ be an arbitrary $2$-chain, then 
\begin{equation}\label{eq:codim2}
\bv(T_1,T_2) = \dfrac{1}{4}-\dfrac{1}{12}\left(\dfrac{d+1-t_2}{d+1-t_1}+\dfrac{t_1}{t_2}\right),
\end{equation}
where $t_i:=|T_i|$ for $i\in\{1,2\}$.
\end{prop}
\begin{proof}
	We use Theorem \ref{thm:main} to compute $\bv(T_1,T_2)$. We apply Equation \eqref{eq:masterformula} to $T_\bullet=(T_1,T_2)$. The set $\mathcal{D}(T_\bullet)$ of all possible spider diagrams is shown in Figure \ref{fig:2}. The leftmost diagram, $\D_1$, consists of trivial spiders. The only contribution of this diagram to the right hand side of \eqref{eq:masterformula} comes from $\Tdcoeff(\D_1)$, since the other statistics are $1$. For the other two diagrams, $\D_2$ and $\D_3$, we have $\Tdcoeff(\D_i)=1/12$, $\sgn(\D_i)=-1$ (notice there is only one leg), $\Bi(\D_i)=1$ for $i\in\{2,3\}$. Finally the corresponding weights for $\D_2,\D_3$ are written in Figure \ref{fig:2}. Formula \eqref{eq:codim2} then follows. 
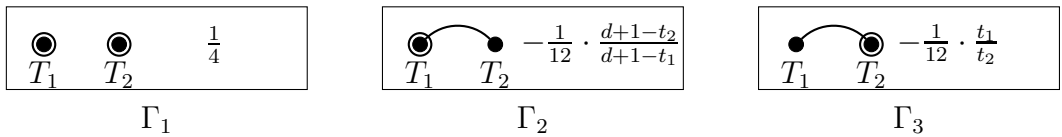
\begin{figure}[ht]
\centering
\begin{tikzpicture}%
\begin{scope}[xshift=0cm,yshift=0cm]
	\draw (.5,.5)--(4.5,.5)--(4.5,-.6)--(.5,-.6)--cycle;
    \node[right] at (3,0) {$\frac{1}{4}$};
    \node at (1,-.4) {$T_1$};
    \node at (2,-.4) {$T_2$};
    \node at (2.5,-1) {$\Gamma_1$};
  \tikzstyle{every node}=[circle, fill, draw,radius=2pt inner sep=1pt,scale=0.5]
    \node (i1) at (1,0) {};
    \node (i2) at (2,0) {};

 	\draw [thick] (1,0) circle [radius=0.15]; 
	\draw [thick] (2,0) circle [radius=0.15];

\end{scope}

\begin{scope}[xshift=5cm,yshift=0cm]
  \draw (.5,.5)--(4.5,.5)--(4.5,-.6)--(.5,-.6)--cycle;
  \node[right] at (2.2,0) {$-\frac{1}{12}\cdot\frac{d+1-t_2}{d+1-t_1}$};
  \node at (1,-.4) {$T_1$};
  \node at (2,-.4) {$T_2$};  
  \node at (2.5,-1) {$\Gamma_2$};
  \tikzstyle{every node}=[circle, fill, draw,radius=2pt inner sep=1pt,scale=0.5]

    \node (i2) at (2,0) {};
    \node (i1) at (1,0) {};

 	\draw [thick] (1,0) circle [radius=0.15]; 
    \draw [thick] (i1) to[out=45, in=135] (i2);

\end{scope}

\begin{scope}[xshift=10cm,yshift=0cm]
	\draw (.5,.5)--(4.5,.5)--(4.5,-.6)--(.5,-.6)--cycle;
	\node[right] at (2.2,0) {$-\frac{1}{12}\cdot\frac{t_1}{t_2}$};
    \node at (1,-.4) {$T_1$};
    \node at (2,-.4) {$T_2$};
    \node at (2.5,-1) {$\Gamma_3$};
  \tikzstyle{every node}=[circle, fill, draw,radius=2pt inner sep=1pt,scale=0.5]
    \node (i1) at (1,0) {};
    \node (i2) at (2,0) {};

 	\draw [thick] (2,0) circle [radius=0.15]; 
    \draw [thick] (i1) to[out=45, in=135] (i2);

\end{scope}

 \end{tikzpicture}
\label{fig:2}
\caption{All spider diagrams on two vertices with the corresponding contribution to Equation \eqref{eq:masterformula}.}
\end{figure}
\end{proof}

Formula \eqref{eq:codim2} (and a similar one for three dimensional cones) was already obtained in \cite{bvalpha} relying on some general formulas in the Berline-Vergne constructions. Since there is no simple closed formula for their construction for unimodular cones of dimension larger than three, we were not able to obtain more formulas in \cite{bvalpha} using the same approach. 
However, by applying Theorem \ref{thm:main}, we obtain in the proposition below a formula for the $\bv$-value of any arbitrary $4$-dimensional braid cone, which could not be obtained with the previously known tools.

\begin{prop}
Let $(T_1,T_2,T_3,T_4)\in\CC^4_{d+1}$ be an arbitrary $4$-chain, then 
\begin{align*}
\bv(T_1,T_2,T_3,T_4) &= \dfrac{1}{16}-\dfrac{1}{48}\left(\dfrac{t_3-t_2}{t_3-t_1}+\dfrac{t_1}{t_2}+\dfrac{t_4-t_3}{t_4-t_2}+\dfrac{t_2-t_1}{t_3-t_1}+\dfrac{d+1-t_4}{d+1-t_3}+\dfrac{t_3-t_2}{t_4-t_2}\right)\\&+\dfrac{1}{144}\left(\dfrac{t_3-t_2}{t_3-t_1}\cdot\dfrac{d+1-t_4}{d+1-t_3}+\dfrac{t_3-t_2}{t_4-t_1}+\dfrac{t_1}{t_2}\cdot\dfrac{d+1-t_4}{d+1-t_3}+\dfrac{t_1}{t_2}\cdot\dfrac{t_3-t_2}{t_4-t_2}\right)\\&+\dfrac{1}{720}\left(\dfrac{t_3-t_2}{t_3-t_1}\cdot\dfrac{t_4-t_3}{t_4-t_1}\cdot\dfrac{d+1-t_4}{d+1-t_1}+3\dfrac{t_1}{t_2}\cdot\dfrac{t_4-t_3}{t_4-t_2}\cdot\dfrac{d+1-t_4}{d+1-t_2}\right)\\&+\dfrac{1}{720}\left(3\dfrac{t_1}{t_2}\cdot\dfrac{t_2-t_1}{t_3-t_1}\cdot\dfrac{d+1-t_4}{d+1-t_3}+\dfrac{t_1}{t_2}\cdot\dfrac{t_2-t_1}{t_4-t_2}\cdot\dfrac{t_3-t_2}{t_4-t_2}\right)
\end{align*}
where $t_i:=|T_i|$ for $i\in\{1,2,3,4\}$.
\end{prop}

\begin{proof}
We use Theorem \ref{thm:main} to compute $\bv(T_1,T_2,T_3,T_4)$. We apply Equation \eqref{eq:masterformula} to $T_\bullet=(T_1,T_2,T_3,T_4)$. The set $\mathcal{D}(T_\bullet)$ of all possible spider diagrams is shown in Figure \ref{fig:4}. Notice that only for two diagrams $\D$ we have a nontrivial $\Bi(\D)$ and in those two cases $\Bi(\D)=3$. Since the coefficients for \eqref{eq:tdexpansion} of odd powers bigger than one are all zero, we do not need to consider spider diagrams with spiders of odd sizes (other than trivial spiders of size one) since $\Tdcoeff$ is zero in that case.
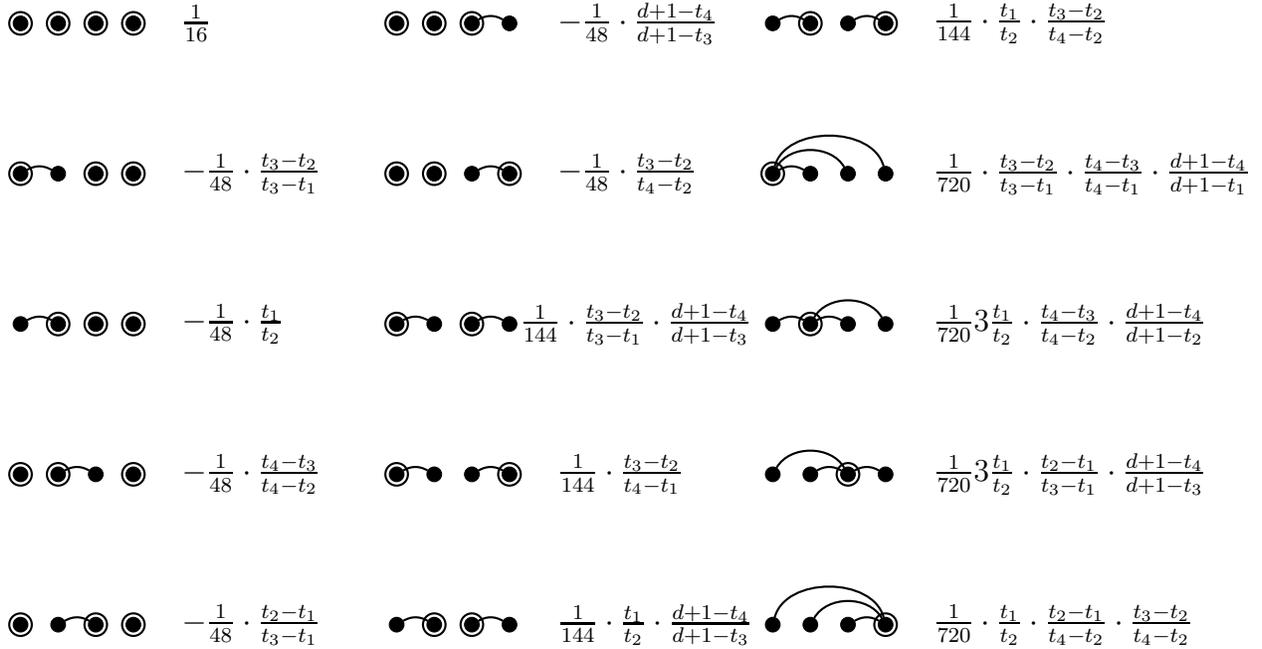
\begin{figure}[ht]
\centering
\begin{tikzpicture}%
\begin{scope}[xshift=0cm,yshift=0cm]

    \node[right] at (3,0) {$\frac{1}{16}$};
  \tikzstyle{every node}=[circle, fill, draw,radius=2pt inner sep=1pt,scale=0.5]
    \node (i1) at (1,0) {};
    \node (i2) at (1.5,0) {};
    \node (i3) at (2,0) {};
    \node (i4) at (2.5,0) {};
    
    \draw [thick] (1,0) circle [radius=0.15];
    \draw [thick] (1.5,0) circle [radius=0.15];
    \draw [thick] (2,0) circle [radius=0.15];
    \draw [thick] (2.5,0) circle [radius=0.15]; 
\end{scope}

\begin{scope}[xshift=0cm,yshift=-2cm]
  \node[right] at (3,0) {$-\frac{1}{48}\cdot\frac{t_3-t_2}{t_3-t_1}$};
  \tikzstyle{every node}=[circle, fill, draw,radius=2pt inner sep=1pt,scale=0.5]

    \node (i2) at (1.5,0) {};
    \node (i1) at (1,0) {};
    \node (i3) at (2,0) {};
    \node (i4) at (2.5,0) {};  
    
    \draw [thick] (1,0) circle [radius=0.15];
    \draw [thick] (2,0) circle [radius=0.15];
    \draw [thick] (2.5,0) circle [radius=0.15];
    \draw [thick] (i1) to[out=30, in=150] (i2);

\end{scope}

\begin{scope}[xshift=0cm,yshift=-4cm]
\node[right] at (3,0) {$-\frac{1}{48}\cdot\frac{t_1}{t_2}$};
  \tikzstyle{every node}=[circle, fill, draw,radius=2pt inner sep=1pt,scale=0.5]
    \node (i1) at (1,0) {};
    \node (i2) at (1.5,0) {};
    \node (i3) at (2,0) {};
    \node (i4) at (2.5,0) {};

    \draw [thick] (1.5,0) circle [radius=0.15];
    \draw [thick] (2,0) circle [radius=0.15];
    \draw [thick] (2.5,0) circle [radius=0.15];

    \draw [thick] (i1) to[out=30, in=150] (i2);

\end{scope}

\begin{scope}[xshift=0cm,yshift=-6cm]
\node[right] at (3,0) {$-\frac{1}{48}\cdot\frac{t_4-t_3}{t_4-t_2}$};
  \tikzstyle{every node}=[circle, fill, draw,radius=2pt inner sep=1pt,scale=0.5]
    \node (i3) at (2,0) {};
    \node (i2) at (1.5,0) {};
    \node (i1) at (1,0) {};
    \node (i4) at (2.5,0) {};      
    \draw [thick] (1,0) circle [radius=0.15];
    \draw [thick] (1.5,0) circle [radius=0.15];
    \draw [thick] (2.5,0) circle [radius=0.15]; 
  
    \draw [thick] (i2) to[out=30, in=150] (i3);

\end{scope}

\begin{scope}[xshift=0cm,yshift=-8cm]
\node[right] at (3,0) {$-\frac{1}{48}\cdot\frac{t_2-t_1}{t_3-t_1}$};
  \tikzstyle{every node}=[circle, fill, draw,radius=2pt inner sep=1pt,scale=0.5]

   \node (i2) at (1.5,0) {};
   \node (i3) at (2,0) {};
   \node (i1) at (1,0) {};  
   \node (i4) at (2.5,0) {};
   
   \draw [thick] (1,0) circle [radius=0.15];
   \draw [thick] (2,0) circle [radius=0.15];
   \draw [thick] (2.5,0) circle [radius=0.15]; 
   \draw [thick] (i2) to[out=30, in=150] (i3);
\end{scope}

\begin{scope}[xshift=5cm,yshift=0cm]
\node[right] at (3,0) {$-\frac{1}{48}\cdot\frac{d+1-t_4}{d+1-t_3}$};
  \tikzstyle{every node}=[circle, fill, draw,radius=2pt inner sep=1pt,scale=0.5]

    \node (i4) at (2.5,0) {};

    \node (i1) at (1,0) {};
    \node (i3) at (2,0) {};
    \node (i2) at (1.5,0) {};

    \draw [thick] (1,0) circle [radius=0.15];
    \draw [thick] (1.5,0) circle [radius=0.15];
    \draw [thick] (2,0) circle [radius=0.15];

    \draw [thick] (i3) to[out=30, in=150] (i4);
\end{scope}

\begin{scope}[xshift=5cm,yshift=-2cm]
\node[right] at (3,0) {$-\frac{1}{48}\cdot\frac{t_3-t_2}{t_4-t_2}$};
  \tikzstyle{every node}=[circle, fill, draw,radius=2pt inner sep=1pt,scale=0.5]

    \node (i3) at (2,0) {};

    \node (i4) at (2.5,0) {};  
    \node (i1) at (1,0) {};
    \node (i2) at (1.5,0) {};

    \draw [thick] (1,0) circle [radius=0.15];
    \draw [thick] (1.5,0) circle [radius=0.15];
    \draw [thick] (2.5,0) circle [radius=0.15]; 
    \draw [thick] (i3) to[out=30, in=150] (i4);
\end{scope}

\begin{scope}[xshift=5cm,yshift=-4cm]
\node[right] at (2.5,0) {$\frac{1}{144}\cdot\frac{t_3-t_2}{t_3-t_1}\cdot\frac{d+1-t_4}{d+1-t_3}$};
  \tikzstyle{every node}=[circle, fill, draw,radius=2pt inner sep=1pt,scale=0.5]

    \node (i2) at (1.5,0) {};
    \node (i4) at (2.5,0) {};

    \node (i1) at (1,0) {};
    \node (i3) at (2,0) {};
    
    \draw [thick] (1,0) circle [radius=0.15];

    \draw [thick] (2,0) circle [radius=0.15];
  
    \draw [thick] (i1) to[out=30, in=150] (i2);
     \draw [thick] (i3) to[out=30, in=150] (i4);
\end{scope}

\begin{scope}[xshift=5cm,yshift=-6cm]
\node[right] at (3,0) {$\frac{1}{144}\cdot\frac{t_3-t_2}{t_4-t_1}$};
  \tikzstyle{every node}=[circle, fill, draw,radius=2pt inner sep=1pt,scale=0.5]

    \node (i2) at (1.5,0) {};
     
    \node (i3) at (2,0) {};

    \node (i1) at (1,0) {};
    \node (i4) at (2.5,0) {}; 
    \draw [thick] (1,0) circle [radius=0.15];
    \draw [thick] (2.5,0) circle [radius=0.15]; 
    \draw [thick] (i1) to[out=30, in=150] (i2);
     \draw [thick] (i3) to[out=30, in=150] (i4);
\end{scope}

\begin{scope}[xshift=5cm,yshift=-8cm]
\node[right] at (3,0) {$\frac{1}{144}\cdot\frac{t_1}{t_2}\cdot\frac{d+1-t_4}{d+1-t_3}$};
  \tikzstyle{every node}=[circle, fill, draw,radius=2pt inner sep=1pt,scale=0.5]
    \node (i1) at (1,0) {};
    \node (i4) at (2.5,0) {};

    \node (i2) at (1.5,0) {};
    \node (i3) at (2,0) {};

    \draw [thick] (1.5,0) circle [radius=0.15];
    \draw [thick] (2,0) circle [radius=0.15];

    \draw [thick] (i1) to[out=30, in=150] (i2);
     \draw [thick] (i3) to[out=30, in=150] (i4);
\end{scope}

\begin{scope}[xshift=10cm,yshift=0cm]
\node[right] at (3,0) {$\frac{1}{144}\cdot\frac{t_1}{t_2}\cdot\frac{t_3-t_2}{t_4-t_2}$};
  \tikzstyle{every node}=[circle, fill, draw,radius=2pt inner sep=1pt,scale=0.5]
    \node (i1) at (1,0) {};
    \node (i3) at (2,0) {}; 
    \node (i2) at (1.5,0) {};
    \node (i4) at (2.5,0) {};

    \draw [thick] (1.5,0) circle [radius=0.15];

    \draw [thick] (2.5,0) circle [radius=0.15]; 
    \draw [thick] (i1) to[out=30, in=150] (i2);
     \draw [thick] (i3) to[out=30, in=150] (i4);
\end{scope}

\begin{scope}[xshift=10cm,yshift=-2cm]
\node[right] at (3,0) {$\frac{1}{720}\cdot\frac{t_3-t_2}{t_3-t_1}\cdot\frac{t_4-t_3}{t_4-t_1}\cdot\frac{d+1-t_4}{d+1-t_1}$};
  \tikzstyle{every node}=[circle, fill, draw,radius=2pt inner sep=1pt,scale=0.5]
    \node (i3) at (2,0) {}; 
    \node (i2) at (1.5,0) {};
    \node (i4) at (2.5,0) {};   
    \node (i1) at (1,0) {};
    
    \draw [thick] (1,0) circle [radius=0.15];

    \draw [thick] (i1) to[out=30, in=150] (i2);
    \draw [thick] (i1) to[out=60, in=120] (i3);
    \draw [thick] (i1) to[out=75, in=105] (i4);
\end{scope}

\begin{scope}[xshift=10cm,yshift=-4cm]
\node[right] at (3,0) {$\frac{1}{720}3\frac{t_1}{t_2}\cdot\frac{t_4-t_3}{t_4-t_2}\cdot\frac{d+1-t_4}{d+1-t_2}$};
  \tikzstyle{every node}=[circle, fill, draw,radius=2pt inner sep=1pt,scale=0.5]
    \node (i1) at (1,0) {};
    \node (i3) at (2,0) {};
    \node (i4) at (2.5,0) {}; 

    \node (i2) at (1.5,0) {};
	\draw [thick] (1.5,0) circle [radius=0.15];
    \draw [thick] (i2) to[out=150, in=30] (i1);
    \draw [thick] (i2) to[out=30, in=150] (i3);
    \draw [thick] (i2) to[out=60, in=120] (i4);
\end{scope}

\begin{scope}[xshift=10cm,yshift=-6cm]
\node[right] at (3,0) {$\frac{1}{720}3\frac{t_1}{t_2}\cdot\frac{t_2-t_1}{t_3-t_1}\cdot\frac{d+1-t_4}{d+1-t_3}$};
  \tikzstyle{every node}=[circle, fill, draw,radius=2pt inner sep=1pt,scale=0.5]
    \node (i1) at (1,0) {};
    \node (i2) at (1.5,0) {};
    \node (i4) at (2.5,0) {}; 
    \node (i3) at (2,0) {};
	\draw [thick] (2,0) circle [radius=0.15];
    \draw [thick] (i3) to[out=150, in=30] (i2);
    \draw [thick] (i3) to[out=30, in=150] (i4);
    \draw [thick] (i1) to[out=60, in=120] (i3);
\end{scope}

\begin{scope}[xshift=10cm,yshift=-8cm]
\node[right] at (3,0) {$\frac{1}{720}\cdot\frac{t_1}{t_2}\cdot\frac{t_2-t_1}{t_4-t_2}\cdot\frac{t_3-t_2}{t_4-t_2}$};
  \tikzstyle{every node}=[circle, fill, draw,radius=2pt inner sep=1pt,scale=0.5]
    \node (i3) at (2,0) {}; 
    \node (i2) at (1.5,0) {};
    \node (i1) at (1,0) {};   
    \node (i4) at (2.5,0) {};
    \draw [thick] (2.5,0) circle [radius=0.15];
    \draw [thick] (i3) to[out=30, in=150] (i4);
    \draw [thick] (i2) to[out=60, in=120] (i4);
    \draw [thick] (i1) to[out=75, in=105] (i4);
\end{scope}
 \end{tikzpicture}
\caption{All spider diagrams on four vertices with the corresponding contribution to Equation \eqref{eq:masterformula}.}
\label{fig:4}
\end{figure}
\end{proof}

\begin{ex}\label{ex:negative}
Using \textsc{sage} \cite{sagemath} we found negative values for $\bv_d$ in four dimensional cones in $\Sigma_d$. The smallest $d$ for which this happens is $d+1=25$, where $\bv_{24}(T_1,T_2,T_3,T_4)=-19/1684800$, for any four chain with $|T_1|=10,|T_2|=12,|T_3|=13, |T_4|=15$. 
\end{ex}

Example \ref{ex:negative} disproves Conjecture \ref{conj:refined}. Furthermore, it also enable us to prove Theorem \ref{thm:noneffective}, which we restate here.

\begin{thm}\label{thm:noneffective-restate}
The Todd class of the permutohedral variety $X_d$ is not effective for $d\geq 24$. That is, there is no way of expressing it as a nonnegative combination of cycles.
\end{thm}
\begin{proof}
It is well known that in the Chow ring of a toric variety arbitrary cycles can be expressed as positive combinations of torus invariant cycles (see Lemma \ref{lem:aux}), so it suffices to show that there is no positive expansion using torus invariant cycles, i.e., that there is no expression of the form \eqref{eq:toddexpression} with all coefficients positive.

By Proposition \ref{prop:symm}, if there is any positive square-free expression for the Todd class of $X_d$, then $\bv(\cdot)$ is positive for all chains $T_\bullet$ in $\CC_{d+1}$, but Example \ref{ex:negative} shows that this is false for $d=24$. Moreover, Remark 3.6 in \cite{bvalpha} implies that there are negative values for all $d\geq 24$. 
\end{proof}

\section{Positivity for linear coefficients}\label{sec:edge} 

As mentioned in the introduction for every lattice polytope $P$ the function $\Lat(tP), t\in \NN$ is a polynomial in $t$ of dimension $d=\dim P$, i.e., $\Lat(tP) = a_0+a_1t^1+a_2t^2+\cdots+a_dt^d,\quad a_i\in\mathbb{Q}.$ This is the \textbf{Ehrhart polynomial} of $P$ and will be denoted $\Lat(P,t)$. We also define $\Lat^i(P):=[t^i]\Lat(P,t)$, the coefficient of $t^i$ in the Ehrhart polynomial.

\subsection{$\bv$ positivity}
In this section we take a different argument to show that $\bv$ values are indeed positive on codimension one cones in the braid fan and thus the main conjecture \ref{conj:main} is true for $\Lat^1$. The arguments in this section are independent of the rest of the paper. We make use of special polytopes called \emph{hypersimplices}.
\begin{defn}
  The \emph{hypersimplex} $\Delta_{k,d+1}$ is defined as 
  \[ \Delta_{k,d+1} = \Perm ( \underbrace{0,\cdots,0}_{d+1-k}, \underbrace{1,\cdots,1}_{k} ).\]
\end{defn}

\begin{prop}\label{prop:hyper}
If $\Lat^1(\Delta_{k,d+1})>0$ for all $1\leq k\leq n$ then $\bv$ is positive on every codimension one cone, thus $\Lat^1(P)>0$ for any generalized permutohedra.
\end{prop}
\begin{proof}
This is a consequence of \cite[Theorem 5.5]{bvalpha}. In the case of an edge the mixed valuation is equal to the valuation itself, the rest of the formula is positive hence the first part follows. The second part is a consequence of the reduction theorem \cite[Theorem 3.5]{bvalpha} which shows how the positivity of $\bv$ for all codimension $k$ cones in $\Sigma_d$ implies positivity of $\Lat^k$ for all generalized permutohedra.
\end{proof}

The following result is standard \cite[Chapter 3, Ex. 62]{ec1}.
\begin{prop}
The Ehrhart polynomial for $\Delta_{k,d+1}$ is given by
\begin{equation}\label{eq:ehrhart}
\Lat(\Delta_{k,d+1}, t) = [z^{kt}]\left(\dfrac{1-z^{t+1}}{1-z}\right)^{d+1}.
\end{equation}
This formula can be turned into the more explicit
\begin{equation}\label{eq:ugly}
\Lat(\Delta_{k,d+1}, t) = \sum_{i=0}^k (-1)^i\binom{d+1}{i}\binom{d+t(k-i)-i}{d}
\end{equation}
\end{prop}

\begin{lem}\label{lem:computation}
For any $k\leq d$, $\Lat^1(\Delta_{k,d+1})>0$.
\end{lem}

\begin{proof}
We are going to keep track of the linear term on each summand on Equation \eqref{eq:ugly}. There are two cases.

For $i=0$ we get $(-1)^0\binom{d+1}{0}\binom{d+tk}{d}$ so
\begin{align*}
[t^1]\binom{d+tk}{d}&=[t^1]\dfrac{(tk+d)\cdots(tk+1)}{d!} =\sum_{i=1}^d\dfrac{k}{i}.
\end{align*}

For $i>0$ we get $(-1)^i\binom{d+1}{i}\binom{d+t(k-i)-i}{d}$ so we first look for the linear term of 
\begin{align*}
\prod_{j=1}^d(t(k-i)-i+j)&=\left(\prod_{j=1}^{i-1}(t(k-i)-i+j)\right)\left(t(k-i)\right)\left(\prod_{j=i+1}^d(t(k-i)-i+j)\right),\\
\end{align*}
which is equal to $(i-1)!(-1)^{i-1}(k-i)(d-i)!$. Now we can compute

\begin{align*}
[t^1](-1)^i\binom{d+1}{i}\binom{d+t(k-i)-i}{d}&=(-1)^i\binom{d+1}{i}\cdot\dfrac{1}{d!} (i-1)!(-1)^{i-1}(k-i)(d-i)!\\
&=-\dfrac{(d+1)(k-i)}{i(d+1-i)}.
\end{align*}

Putting these equations together we get that the linear term in \eqref{eq:ugly} is
\begin{align*}
[t^1]\sum_{i=0}^k (-1)^i\binom{d+1}{i}\binom{d+t(k-i)-i}{d}&=\sum_{i=0}^k [t^1]\left((-1)^i\binom{d+1}{i}\binom{d+t(k-i)-i}{d}\right),\\
&=\sum_{i=1}^n\dfrac{k}{i} - \sum_{i=1}^k\dfrac{(d+1)(k-i)}{i(d+1-i)},\\
&=\sum_{i=1}^k \left(\dfrac{k}{i}-\dfrac{(d+1)(k-i)}{i(d+1-i)}\right) + \sum_{i=k+1}^n \dfrac{k}{i},\\
&>0,
\end{align*}
because each parenthesis is positive since $\dfrac{k}{k-i}\geq \dfrac{d+1}{d+1-i}$ for $d+1>k$. 
\end{proof}

\begin{thm}\label{thm:linear}
Conjecture \ref{conj:main} is true for the linear terms. More precisely, $\Lat^1(P)>0$ for every lattice generalized permutohedron $P$.
\end{thm}
\begin{proof}
It follows from Proposition \ref{prop:hyper} and Lemma \ref{lem:computation}.
\end{proof}
In \cite{kat} the authors use their results in Minkowski linear functionals to give an alternative proof of Theorem \ref{thm:linear}. 


\appendix 
\section{Algebraic cycles in toric varieties} 

We include a sketch of the following lemma, since it is a crucial reduction in the proof of Theorem \ref{thm:noneffective-restate} and we could not find a reference in the literature.
\begin{lem}\label{lem:aux}
Let $X$ be smooth projective toric variety of dimension $n$ over an algebraically closed field $\k$. Then arbitrary algebraic cycles of $X$ are \emph{positive} combination of the torus invariant cycles.
\end{lem}

\begin{proof}[Sketch of proof:]
 Let $T\subset X$ be its dense torus and fix an isomorphism $T\cong T_1\times\cdots\times T_n$, where each $T_i$ is isomorphic to $\k^*$.
 
 Let $Z$ be a general cycle. We can assume that $Z$ is an irreducible subvariety. Taking the \emph{flat limit} (See \cite[Section II.3.4]{schemes}) of $Z$ over $T_1\cong \k^*$ we obtain an subvariety $Z_1$ whose associated algebraic cycle is effective and rationally equivalent to $Z$, and furthermore each irreducible component of $Z_1$ is $T_1$-invariant. For every irreducible component of $Z_1$ we now take the flat limit over $T_2$, taking the union we obtain an effective cycle $Z_2$ rationally equivalent to $Z_1$ that is $T_1$- and $T_2$-invariant. Continuing in this way after $n$ iterations we get at an effective cycle $Z_n$ rationally equivalent to each $Z_i, i<n$, and to the original $Z$, that is $T_i$-invariant for every $i$. Then $Z_n$ is rationally equivalent to $Z$ and $T$-invariant as we wanted to show.
\end{proof}

The intuition is that for each tori $T_i\cong\k^*$ we have an action on any subvariety $Z\subset X$, so there exist subvarieties $t\cdot Z\subset X$ for any $t\in\k^*$ (which are all isomorphic to $Z$) and we take the limit as $t$ approaches $0$ to obtain a rationally equivalent subvariety that is now $\k^*$- invariant.

\begin{ex}
Let $X=\mathbb{P}^3_\k$ be a toric variety with torus $T=\{(t_1:t_2:t_3:1):(t_1,t_2,t_3)\in\left(\k^*\right)^3\}$ acting coordinate-wise. Consider $Z=V(xy-z^2-w^2,xw-yz)\subset \mathbb{P}^3_\k$ to be the surface given by the zero locus of system of equations 
\begin{align*}
xy&=z^2+w^2,\\
xw&=yz.
\end{align*} 
One can check that $Z$ is irreducible. We compute $Z_1,Z_2$ and $Z_3$ as in the proof of Lemma \ref{lem:aux}.

\begin{enumerate}

	\item We have that $T_1:=\{(t_1:1:1:1):t_1\in\k^*\}$ acts by scaling the first coordinate. For a fixed nonzero scalar $t\in\k^*$, the subvariety $(t:1:1:1)\cdot Z$ is equal to $V(xy-tz^2-tw^2, xw-tyz)$. Taking the flat limit as $t\to 0$, we obtain \[Z_1:=V(xy,xw,y^2z-z^2w-w^3).\] This subvariety decomposes as $Z_1=U\cup W$, where $U=V(y,w)$ and  $W=V(x,y^2z-z^2w-w^3)$. Notice that both components are $T_1$-invariant. In the Chow ring we get \[[Z_1]=[U]+[W].\]
The cycle $[Z_1]$ is rationally equivalent to $[Z]$.
\item Next we have $T_2:=\{(1:t_2:1:1):t_2\in\k^*\}$ and it acts by scaling the second coordinate. 
	For a fixed nonzero scalar $t\in\k^*$, the subvariety $(1:t:1:1)\cdot W$ is equal to $V(x,y^2z-t^2z^2w-t^2w^3)$. Taking the flat limit as $t\to 0$ we obtain $V(x,y^2z)$. This subvariety has two components, $U'=V(x,y)$ and $U''=V(x,z)$ where $U'$ has multiplicity two. Hence, in the Chow ring we have 
	\[ [W]=[V(x,y^2z)]=2[U']+[U''].\]
The component $U$ is already $T_2$-invariant, so it is equal to its flat limit over $T_2$. Collecting terms from both components we get $Z_2:=V(x,y^2z)\cup V(y,w)$ which can be represented in the Chow ring as
\[[Z_2]=[U]+2[U']+[U''].\]
The algebraic cycle $[Z_2]$ is rationally equivalent to $[Z_1]$ and hence to $[Z]$.

\item Finally we have $T_3:=\{(1:1:t_3:1):t_3\in\k^*\}$ acting by scaling the third coordinate. Since every irreducible component of $Z_2$ is already $T_3$-invariant, each one is equal to its flat limit over $T_3$, so nothing changes in this step and $Z_3=Z_2$ so in the Chow ring \[[Z_3]:=[U]+2[U']+[U''].\]
\end{enumerate}
We obtained the expression $[V(y,w)]+2[V(x,y)]+[V(x,z)]$ in the Chow ring $A(\mathbb{P}^3_\k)$ which is a positive combination of three torus invariant cycles and it is rationally equivalent to $[Z]$ as we wanted.
\end{ex}

\bibliography{biblio}

\def\cprime{$'$}
\begin{thebibliography}{10}

\bibitem{bvtodd}
Nicole Berline and Mich{\`e}le Vergne.
\newblock The equivariant {T}odd genus of a complete toric variety, with
  {D}anilov condition.
\newblock {\em J. Algebra}, 313(1):28--39, 2007.

\bibitem{bvoriginal}
Nicole Berline and Mich\`ele Vergne.
\newblock Local {E}uler-{M}aclaurin formula for polytopes.
\newblock {\em Mosc. Math. J.}, 7(3):355--386, 573, 2007.

\bibitem{bvalpha}
Federico Castillo and Fu~Liu.
\newblock Berline-{V}ergne valuation and generalized permutohedra.
\newblock {\em Discrete Comput. Geom.}, 60(4):885--908, 2018.

\bibitem{nested}
Federico Castillo and Fu~Liu.
\newblock Deformation cones of nested braid fans.
\newblock {\em arxiv:1710.01899}, 2018.

\bibitem{nill}
Federico Castillo, Fu~Liu, Benjamin Nill, and Andreas Paffenholz.
\newblock Smooth polytopes with negative {E}hrhart coefficients.
\newblock {\em J. Combin. Theory Ser. A}, 160:316--331, 2018.

\bibitem{danilov}
V.I. Danilov.
\newblock The geometry of toric varieties.
\newblock {\em Russian Math. Surveys}, 33(2):97--154, 1978.

\bibitem{deloera}
Jes{\'u}s~A. De~Loera, David~C. Haws, and Matthias K{\"o}ppe.
\newblock Ehrhart polynomials of matroid polytopes and polymatroids.
\newblock {\em Discrete Comput. Geom.}, 42(4):670--702, 2009.

\bibitem{ehrhart}
Eug{\`e}ne Ehrhart.
\newblock Sur les poly\`edres rationnels homoth\'etiques \`a {$n$}\ dimensions.
\newblock {\em C. R. Acad. Sci. Paris}, 254:616--618, 1962.

\bibitem{schemes}
David Eisenbud and Joe Harris.
\newblock {\em The geometry of schemes}, volume 197 of {\em Graduate Texts in
  Mathematics}.
\newblock Springer-Verlag, New York, 2000.

\bibitem{ewald}
G{\"u}nter Ewald.
\newblock {\em Combinatorial convexity and algebraic geometry}, volume 168 of
  {\em Graduate Texts in Mathematics}.
\newblock Springer-Verlag, New York, 1996.

\bibitem{hypersimplices}
Luis Ferroni.
\newblock Uniform matroids are ehrhart positive.
\newblock {\em arXiv:1911.10146}, 2019.

\bibitem{fulton}
William Fulton.
\newblock {\em Introduction to toric varieties}, volume 131 of {\em Annals of
  Mathematics Studies}.
\newblock Princeton University Press, Princeton, NJ, 1993.
\newblock The William H. Roever Lectures in Geometry.

\bibitem{garpomm}
Stavros Garoufalidis and James Pommersheim.
\newblock Sum-integral interpolators and the {E}uler-{M}aclaurin formula for
  polytopes.
\newblock {\em Trans. Amer. Math. Soc.}, 364(6):2933--2958, 2012.

\bibitem{kat}
Katharina Jochemko and Mohan Ravichandran.
\newblock Generalized permutahedra: Minkowski linear functionals and ehrhart
  positivity.
\newblock {\em arXiv:1909.08448}, 2019.

\bibitem{positivesurvey}
Fu~Liu.
\newblock On positivity of {E}hrhart polynomials.
\newblock In {\em Recent trends in algebraic combinatorics}, volume~16 of {\em
  Assoc. Women Math. Ser.}, pages 189--237. Springer, Cham, 2019.

\bibitem{mcmullenext}
Peter McMullen.
\newblock Valuations and dissections.
\newblock In {\em Handbook of convex geometry, {V}ol. {A}, {B}}, pages
  933--988. North-Holland, Amsterdam, 1993.

\bibitem{toddclass}
James Pommersheim and Hugh Thomas.
\newblock Cycles representing the {T}odd class of a toric variety.
\newblock {\em J. Amer. Math. Soc.}, 17(4):983--994, 2004.

\bibitem{ec2}
Richard~P. Stanley.
\newblock {\em Enumerative combinatorics. {V}ol. 2}, volume~62 of {\em
  Cambridge Studies in Advanced Mathematics}.
\newblock Cambridge University Press, Cambridge, 1999.
\newblock With a foreword by Gian-Carlo Rota and appendix 1 by Sergey Fomin.

\bibitem{ec1}
Richard~P. Stanley.
\newblock {\em Enumerative combinatorics. {V}olume 1}, volume~49 of {\em
  Cambridge Studies in Advanced Mathematics}.
\newblock Cambridge University Press, Cambridge, second edition, 2012.

\bibitem{sagemath}
{The Sage Developers}.
\newblock {\em {S}ageMath, the {S}age {M}athematics {S}oftware {S}ystem
  ({V}ersion x.y.z)}, YYYY.
\newblock {\tt https://www.sagemath.org}.

\end{thebibliography}
\bibliographystyle{plain}

\end{document}